# Distributed Optimal Power Flow

HyungSeon Oh[1*]

[1] Department of Electrical and Computer Engineering, United States Naval Academy, Annapolis, Maryland, United States of America
* Email: hoh@usna.edu

# Abstract

**Objective:** The objectives of this paper are to 1) construct a new network model compatible with distributed computation, 2) construct the full optimal power flow (OPF) in a distributed fashion so that an effective, non-inferior solution can be found, and 3) develop a scalable algorithm that guarantees the convergence to a local minimum.

**Existing challenges:** Due to the nonconvexity of the problem, the search for a solution to OPF problems is not scalable, which makes the OPF highly limited for the system operation of large-scale real-world power grids—"the curse of dimensionality". The recent attempts at distributed computation aim for a scalable and efficient algorithm by reducing the computational cost per iteration in exchange of increased communication costs.

**Motivation:** A new network model allows for efficient computation without increasing communication costs. With the network model, recent advancements in distributed computation make it possible to develop an efficient and scalable algorithm suitable for large-scale OPF optimizations.

**Methods:** We propose a new network model in which all nodes are directly connected to the center node to keep the communication costs manageable. Based on the network model, we suggest a nodal distributed algorithm and direct communication to all nodes through the center node. We demonstrate that the suggested algorithm converges to a local minimum rather than a point, satisfying the first optimality condition.

**Results:** The proposed algorithm identifies solutions to OPF problems in various IEEE model systems. The solutions are identical to those using a centrally optimized and heuristic approach. The computation time at each node does not depend on the system size, and $N_{iter}$ does not increase significantly with the system size.

**Conclusion:** Our proposed network model is a star network for maintaining the shortest node-to-node distances to allow a linear information exchange. The proposed algorithm guarantees the convergence to a local minimum rather than a maximum or a saddle point, and it maintains computational efficiency for a large-scale OPF, scalable algorithm.

# Introduction

In modern societies, demand for electricity is expected to be satisfied continuously via controllable generation technologies. An event is a situation in which the demand is not fulfilled. Ten-in-one, a widely used reliability criterion for events, means that an event should occur just once in a 10-year span. To meet this standard, system operators schedule the generation portfolio and the grid systems in advance. For example, a day-ahead unit commitment determines the 24-hourly dispatches, along with unit commitment decisions, to meet varying hourly demands. For each hour, the demand profiles are assumed to be constant, which defines the process's steady-state operation. In the absence of an unexpected disturbance, stochastic hourly demand is the unique source of uncertainty in traditional power system operation. Over the last decade or so, renewable energy resources and smart grid technologies have been integrated into systems to improve energy efficiency and reduce greenhouse gas emissions. This integration has introduced uncertainty into the operation of power systems, presenting a new challenge. If high-precision forecasting could be introduced to estimate future energy resources and control demand, existing operations' tools would remain useful, assuming they could be integrated into the expected effective demand (≡ demand – expected demand reduction – expected renewable energy resources). Unfortunately, even though the precision of forecasting tools has improved, the errors in their long-term forecasts, for a day ahead, for example, are not yet sufficiently small for reliable operation. Frequent decisions are a potential way to accommodate the uncertainty. For example, a day-ahead 24-hour unit commitment (UC) decision is made once in a daily cycle. If the errors in 2-hour ahead forecasts are small enough, then the UC decision with a forecast every 2 hours would still be a reliable tool for the power system's operation. The computational capability to support such decisions plays a key role in this process.

Optimal power flow (OPF) is a backbone in the steady-state operation of power systems. The characteristics of OPF are highly nonlinear and nonconvex. The computational complexity associated with these characteristics makes power flow (PF) analysis a non-deterministic polynomial-time (*NP*)-hard problem [1]. In most operational practices, a linear approximation of OPF, namely, direct current (DC) OPF, is pursued. Although easy to solve, DC OPF does not address voltage problems, losses, and the dispatch of reactive power generation. Due to these issues, DC OPF may not be feasible. To address the problem correctly, it is ideal to aim for a nonlinear and nonconvex OPF. In addition to the nonconvex nature of the full OPF, uncertainties increase the number of variables in traditional, central decision-making processes. Therefore, frequent but short-term decisions concerning large-scale power systems can be challenging. With the recent advancements in hardware in multi-core machines, distributed computation becomes an attractive approach for enhancing computational efficiency. An exemplary area in power system analysis concerns the use of distributed computation for OPF. Motivated readers can find information related to distributed approaches to solving OPF problems [2,3]. Within distributed computation, the alternating direction method of multipliers (ADMM) has gained popularity due to its straightforward implementation and its provable convergence (if the original problem is convex). For a linear DC OPF, the ADMM approach can result in successful convergence to a global solution [4–6]. A full OPF problem on a radial network can be relaxed to a convex semidefinite programming (SDP) problem, and the relaxed problem is exact if optimal power injections lie in a region where the voltage upper bounds do not bind [7,8]. Several studies have used the ADMM approach to solve OPF problems in a radial network [9–11]. OPF problems usually involve the operation of mesh transmission systems, but an SDP solution for a mesh network may not be physically meaningful [2].

A nodal OPF would be the most intuitive approach to extending a central OPF in a distributed fashion. In the nodal OPF [12], the information exchange among the nearest neighbors leads to high communication costs. The maximum node-to-node distance (also termed the path length (PL) [13]) plays a key role in the communication costs. The convergence tends to be very slow due to the contaminated information received from local decisions during communications. To the best of our knowledge, there has been no report of the successful convergence of this approach for any mesh networks.

The PL can be reduced when a clustering approach is undertaken for a mesh network via the partitioning of a system into multiple subsystems [14–19]. Two adjacent subsystems share some nodes and branches between the nodes; thus, the PL is small, keeping the communication costs manageable. The primary and the dual variables at the shared nodes and branches are constrained equally. This approach can be efficient if the shared nodes adequately represent the other nodes in the same subsystem. Several studies have proposed an efficient algorithm for partitioning a system so that the ADMM converges to a solution [16,17,19,20]. In contrast to a study by Sun, Phan, and Ghosh [12], the flow constraints can be integrated for the lines of intra-subsystems. In two studies, Erseghe [14,15] integrated the flow limits of the lines of inter-subsystems by redefining the subsystems to overlap the lines. Guo, Hug, and Tonguz [16] report that the inclusion of the limit still yields a solution. However, the approach has several shortcomings that contribute to computational inefficiency: 1) the low quality of the solution, 2) the need for a warm starting point for convergence, and 3) the communication costs. In addition, the convergence behavior is not reported, so it is not possible to discuss the computational efficiency. The solutions presented in several studies [14–17] are low-voltage, inferior solutions due to increased losses. Engelmann et al. [18] added a significantly large term regarding reactive power injections that affects the optimality conditions and, as a result, the distributed problem is different from the original one. The necessity of a solution for the nonconvex PF as a starting point increases the computational costs. Even though the PL decreases in comparison to the nodal OPF in the study by Sun, Phan, and Ghosh [12], the communication costs increase the computational costs due to the tradeoff discussed by Guo, Hug, and Tonguz [16]. These shortcomings make the benefit of distributed computation questionable. The distributed SDP approach by Madani, Kalbat, and Lavaei [19] yields the global (and therefore identical) solution to the central SDP approximation to the OPF problem because the problem is convex, but the solution may not be physically feasible. In addition to the computational inefficiency, there is no approach that can theoretically yield an optimal balance between the computational cost for one subsystem and the communication costs of the subsystems. Another study by Guo, Hug, and Tonguz [20] proposes a heuristic approach for selecting subsystems, but it does not yield a unique choice because its initialization is based on the local solution for the nonconvex OPF. In several studies [14,17,18,20], a positive correlation was observed between the system size and the number of nodes in the largest subsystems (see Fig 1). The PL depends linearly on the number of subsystems *n*, which depends inversely on the size of the largest subsystem ($N_{sub}$), i.e., $n \geq Nb/N_{sub}$, where the equality holds when the sizes of the subsystems are uniform. Fig 1 indicates a positive correlation between $N_{sub}$ and *Nb*. The dotted line indicates that the best linear fit for the relationship is $N_{sub} \propto Nb^{0.88}$. If the optimization problem is solved by SDP, the computation time is found in $\vartheta\left(N_{sub}^{3}\right)$ [21] that yields a computation time that is proportional to $\vartheta\left(Nb^{2.64}\right)$. In addition, the approaches require significant communication costs as well. Therefore, the overall computation cost is much higher than that of the central optimization. In addition to the computational efficiency, it is not guaranteed that the voltages at the boundary of each subsystem represent the voltages of other nodes inside the same subsystem correctly. If they do not, convergences may not be observed because the information exchange is limited to the boundary buses. The slow convergences and/or non-monotonous convergences reported in previous studies [14–17] indicate the insufficient representativeness of the boundary buses. A relatively fast convergence is reported by Engelmann et al. [18] in exchange of per-step communication costs by sharing the sensitivities in addition to the local primal variables. The increase in the communication cost is found in $\vartheta\left(\sum_i n_i^2\right)$ where $n_i$ is the dimension of the local gradient at the $i^{th}$ group. Although the progress at each iteration is faster than those of other ADMM approaches, the communication cost itself is much higher than the total computation costs of the nonconvex heuristic solvers or of the SDP solvers. As a result, although these studies are worth exploring, we conclude that the aggregation approach is not practical in terms of computational efficiency and the inferior quality of the solutions for a scalable algorithm due to the tradeoff issue, nonconvexity, and the modeling problem. A new network model for OPF is necessary for the distributed computation.

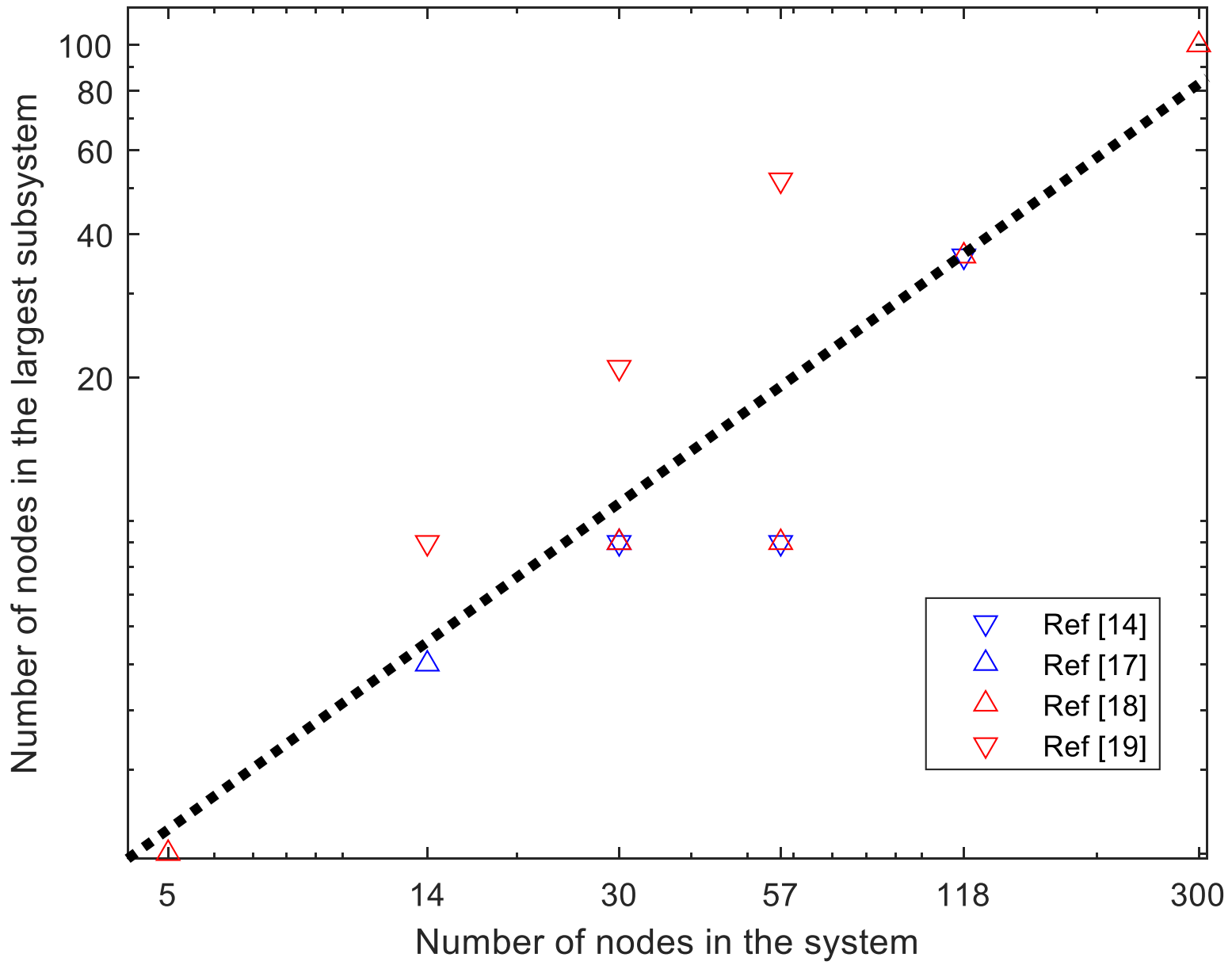

Fig 1. The largest number of nodes in the partitioned subsystems of various power system test cases reported in the literature [14, 17-19].

The contributions of this paper are 1) a new network model that yields direct communication among nodes regardless of the system size, 2) a distributed, fast, and efficient algorithm to solve highly nonlinear and nonconvex

OPF problems, and 3) a scalable algorithm that guarantees convergence to a local minimizer. The paper is organized as follows: the theory section proposes a new network model designed for distributed computation and presents an algorithm to solve a full OPF; the next section describes the details of the implementation of the proposed algorithm; the results and discussion present the results for the OPF and comparisons with those from other studies; the following section provides conclusions and research directions for further improving the computational efficiency; and the appendices sketch the proofs of the ranks of the matrices associated with OPF problems and of the convergence.

# Theory

## Proposed network model and algorithm

We propose a new network model for a nodal OPF to keep communication costs manageable regardless of the system size. For this purpose, the desired properties of the model are as follows:

1) The model must be compatible with PF studies for which Kirchhoff's laws and voltage magnitudes are well defined.
2) Each node is a short distance away from the rest of the nodes to minimize the communication costs.
3) The voltages at a node and those at the rest of the nodes are linearly related.

In the power flow studies, such as PF, OPF, state estimation, and probabilistic PF, voltages are the variables. The constraints in the studies consist of the power flows and injections, as well as the voltage magnitudes in terms of voltages. In the Cartesian coordinate system, the power flows and injections are quadratic in voltage. For example, the power flow over the line connecting Nodes $i$ and $j$ at $i$ is $\boldsymbol{f}_{i-j}^{i} = \boldsymbol{v}_i \boldsymbol{i}_{i-j}^{*} = v^T \boldsymbol{F}_{i-j}^{i} v$, and the power injection at Node $i$ is $\boldsymbol{g}_i = p_i + \mathbf{j} q_i = v^T \boldsymbol{S}_i v$ where the quantities sandwiched between voltages are in $2Nb$-by-$2Nb$; $\boldsymbol{F}_{i-j}^{i} \triangleq \boldsymbol{J} e_i e_{i-j}^T \boldsymbol{Y}_{br}^{*} \boldsymbol{J}^H$ where $e_{i-j}$ is a vector with the cardinality of $2Nl$ of which the element corresponding to the flow over $i$-$j$ is 1, and all other elements are zeros, $\boldsymbol{S}_i \triangleq \boldsymbol{J} e_i e_i^T \boldsymbol{Y}_{bus}^* \boldsymbol{J}^H$ ; and the superscript $H$ is the conjugate transpose. The matrices $\boldsymbol{F}_{i-j}^i$ and $\boldsymbol{S}_i$ have two nonzero rows at $i$ and $Nb+i$ rows due to $\boldsymbol{J}e_i$.

**Claim 1**: The matrices associated with power flows and power injections are all of rank 4.

**Claim 2**: The matrices associated with the squares of the voltage magnitudes have rank 2.

See Appendix A for the proof of the claims.

For a real-valued symmetric rank-4 matrix $M_j$, a real-valued eigen pair $\lambda$ and $u$ exist such that $M_j = \sum_{k=1}^{4} \lambda_{jk} u_{jk} u_{jk}^T = \underbrace{\left[ \sqrt{|\lambda_{j1}|} u_{j1} \; \sqrt{|\lambda_{j2}|} u_{j2} \; \sqrt{|\lambda_{j3}|} u_{j3} \; \sqrt{|\lambda_{j4}|} u_{j4} \right]}_{\Phi_j} II_j \Phi_j^T$. Because $M_j$ is a rank 4 matrix, the number of columns in $\Phi_j$ and the number of nonzero diagonal elements in $II_j$ are both 4. Here, the terms $\Phi_{\breve{S}_j}, \Phi_{\hat{S}_j}, \Phi_{F_{j-l}^j}$, and $\Phi_{\hat{F}_{j-l}^j}$ : $\Phi_{\breve{S}_j}$ refer to the matrices associated with the real power injection at Node $j$ $\left( \breve{S}_j \triangleq \mathrm{Re}\left(\boldsymbol{S}_j\right) = \Phi_{\breve{S}_j} II_j^{\breve{S}_j} \Phi_{\breve{S}_j}^T \right)$; $\Phi_{\hat{S}_j}$ is associated with reactive power injection at Node $j$ $\left( \bar{S}_j \triangleq \mathrm{Im}\left(\boldsymbol{S}_j\right) = \Phi_{\bar{S}_j} II_j^{\bar{S}_j} \Phi_{\bar{S}_j}^T \right)$; $\Phi_{F_{j-l}^j}$ is associated with real power flow over the line $j$-$l$ at the side of Node $j$ $\left( F_{j-l}^j = \Phi_{F_{j-l}^j} II_{j-l}^{F_{j-l}^j} \Phi_{F_{j-l}^j}^T \right)$; and $\Phi_{\hat{F}_{j-l}^j}$ is associated with the reactive power flow over the line $j$-$l$ at the side of Node $j$ $\left( \bar{F}_{j-l}^j = \Phi_{\bar{F}_{j-l}^j} II_{j-l}^{\bar{F}_{j-l}^j} \Phi_{\bar{F}_{j-l}^j}^T \right)$. Note that the eigenpairs are dependent solely on the system parameters (not on voltages). With the eigenpairs, the quantities of interest in the PF are the real power injection at Node $j$ $p_j$, the reactive power injection at Node $j$ $q_j$, the real power flow over a line $j$-$l$ at Node $j$ $f_{j-l}^j$, the reactive power flow over a line $j$-$l$ at Node $j$ $\bar{f}_{j-l}^j$, and the squared voltage's magnitude at Node $j$ $E_j$. There are quadratic relationships among them:

$$p_j \left(= g_j - d_j\right) = \alpha_j^T II_j^{S_j} \alpha_j ; \; q_j \left(= \bar{g}_j - \bar{d}_j\right) = \beta_j^T II_j^{\bar{S}_j} \beta_j ; \; f_{j-l}^j = \gamma_{j-l,j}^T II_{j-l,j}^{F_{j-l}^j} \gamma_{j-l,j} ; \; \bar{f}_{j-l}^j = \delta_{j-l,j}^T II_{j-l,j}^{\bar{F}_{j-l}^j} \delta_{j-l,j} ; \; E_j = \omega_j^T \omega_j \quad (1)$$

Eq. (1) expresses Kirchhoff's laws and the voltage magnitudes. The nodal variables $\alpha_j$, $\beta_j$, $\gamma_{j-l,j}$, $\delta_{j-l,j}$, and $\omega_j$ are defined as follows: the $\alpha_j$ are the voltages projected onto the eigenspaces spanned by the real power injection, and the $\beta_j$ cover the reactive power injection, the $\gamma_{j-l,j}$ cover the real power flow over $j$-$l$ at Node $j$, the $\delta_{j-l,j}$ cover the reactive power flow over $j$-$l$ at Node $j$, and the $\omega_j$ cover the voltage magnitudes at the $j^{th}$ node in the same manner. Although there is one each of $\alpha_j$, $\beta_j$, and $\omega_j$ at each node, the $nl_j$ of $\gamma_{j-l,j}$ and $\delta_{j-l,j}$ are defined such that $nl_j$ represents the number of branches that are directly connected to Node $j$ (i.e., $l$ = 1, 2, …, $nl_j$). It is worth noting that the nodal variables are linearly associated with the voltages (the definitions of the nodal variables in Eq. (1)), which indicates that PL = 1. Therefore, the nodal variables satisfy all the desired properties of the new network model.

Next, we let $nl_j$ be the number of lines connected to the $j^{\text{th}}$ node. Then, the dimensions of α, β, γ, δ, and ω are 4, 4, $4nl_j$, $4nl_j$, and 2, respectively. Note that their dimensions do not depend on the system size. Let $\mu_j$ be a nodal variable that integrates local generation variables as follows: $\mu_j = \left[ x_j^T \ {}_Q x_j^T \ \mu_{end} \right]^T$; $x_j^{8nl_j+10} = \left[ \alpha_j^T \ \beta_j^T \ \gamma_{jl}^T \ \delta_{jl}^T \ \omega_j^T \right]^T$; ${}_Q x_j^{2nl_j+2ng_j} = \left[ \left( f_{j-l}^{j} \right)^T \left( \bar{f}_{j-l}^{j} \right)^T \ p_j^T \ q_j^T \right]^T$ where $\mu_{end}$ = 1. The generalized PF formulation is listed in previous work [22]. Note that the cardinality of the variable is fixed—independent of the system size. Therefore, the new network model yields: 1) a fixed number of variables, and 2) a PL of 1.

# Proposed network model: a star grid with two channels

In the definition of the local variable $x_j$, it is noticeable that there are three variable types: 1) power injection and flow variables α, β, γ, and δ, 2) the voltage magnitude variable ω, and 3) quadratic variables, ${}_Q x_j$ comprising power flows $\left( f_{j-k}^{j}, \bar{f}_{j-k}^{j} \right)$ and power generations ($g_j$). The first two types of variables are linear, and the third type is quadratic, in terms of the voltages. Once $x_j$ (α, β, γ, and δ) is determined, it is straightforward to find the power generations at Node $j$ with Eq. (1) (i.e., $x_j$ defines the feasible region of $\mu_j$). The voltage $v_p$ can be reconstructed from α, β, γ, and δ. The communication path used to collect their values is termed the power channel. Here, ω is directly associated with the voltage magnitudes, and the values are reported through another communication path termed the voltage channel to update voltages $v_m$. Although the voltages $v_p$ and $v_m$ should be identical, conditions are relaxed so that they can be different.

Fig 2 illustrates an example of the traditional network model (left) and the proposed network model (right) for a modified 4-bus system that has branches connecting 1–2, 1–3, 1–4, 2–4, and 3–4. In the proposed model, all the nodes are connected to the center node (black dot at the center) through the power channel (red lines) and the voltage channel (green lines), and all the nodes are of distance 1 from the center (PL = 1). The network is a complete bipartite graph of order 2 with a maximum-diameter star network. All communications between the variables take place in local nodes ($\alpha$, $\beta$, $\gamma$, $\delta$, and $\omega$), and the center node is linear in terms of voltages.

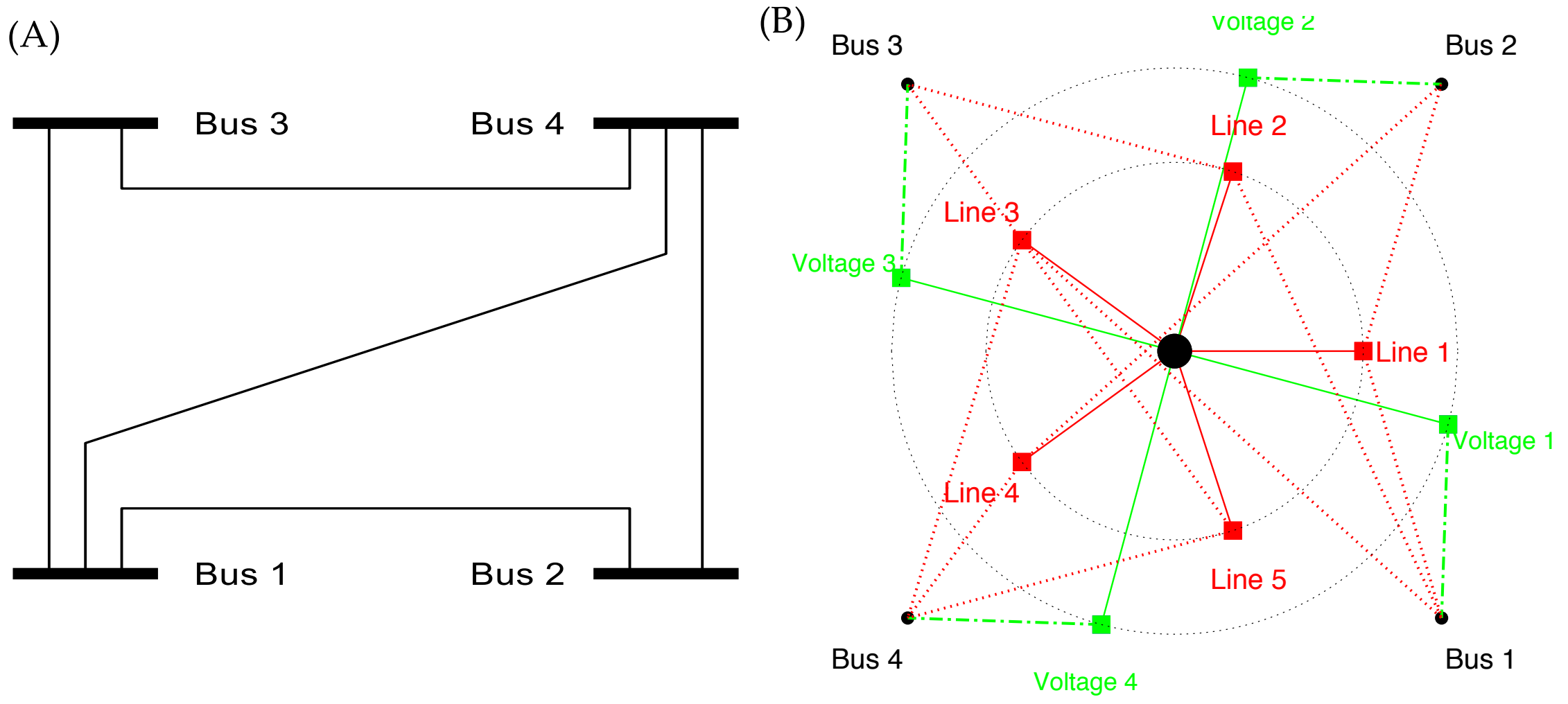

Fig 2. A star and linear network for a modified IEEE-4 bus system. The traditional network model is shown on the left (A) and the proposed network model is presented on the right (B). Red lines indicate the power channel, and green lines are the voltage channel.

# Nodal OPF with nodal variables and its convex relaxation

Even though we describe the OPF problem in this paper, the same framework is applied to any other problem, including PF, OPF, state estimation, and probabilistic PF problems, because they have unified formulas [22]. The central OPF formulation as a nonconvex, quadratically constrained quadratic problem is given in Eq. (3), its distributed optimization problem (nodal OPF) at each node $I$ is given in Eq. (4A), and the nodal OPF is formulated in terms of an ADMM algorithm in Eq. (4B).

$$\min_{v,g} g^T diag(a) g + 2b^T g, \quad g = \left[ g^T \ \ \bar{g}^T \right]^T \tag{3}$$

$$\text{s.t.} \quad \begin{array}{ll}
1. & v^T S_j v - L_g^T g + d_j = 0, v^T \bar{S}_j v - L_g^T \bar{g} + \bar{d}_j = 0 \quad \forall j \in B \\
2. & v^T F_{i-j}^i v - f_{i-j}^i = 0, v^T \bar{F}_{i-j}^i v - \bar{f}_{i-j}^i = 0 \quad \forall j \in Br\{i\} \\
3. & f^T \left( e_l e_l^T + e_{l+Nl} e_{l+Nl}^T \right) f - cap_l^2 \le 0, f = \left[ f^T \ \bar{f}^T \right]^T \quad \forall l \in Br \\
4. & g^T \left( e_m e_m^T \right) g - \left( p_m^{\max} + p_m^{\min} \right) e_m^T g + p_m^{\max} \cdot p_m^{\min} \le 0 \quad \forall m \in G \\
5. & g^T \left( e_{m+Ng} e_{m+Ng}^T \right) g - \left( q_m^{\max} + q_m^{\min} \right) e_{m+Ng}^T g + q_m^{\max} \cdot q_m^{\min} \le 0 \\
6. & v^T \left( e_j e_j^T + e_{j+Nb} e_{j+Nb}^T \right) v - E_j^{\max} \le 0 \quad \forall j \in B \\
7. & -v^T \left( e_j e_j^T + e_{j+Nb} e_{j+Nb}^T \right) v + E_j^{\min} \le 0 \quad \forall j \in B
\end{array}$$

$$\min_{\mu} L(\mu) = \min_{\mu} \sum_i L_i(\mu_i) = \min_{\mu} \sum_i \mu_i^T M_{ob}^i \mu_i \tag{4A}$$

where $M_{ob}^i = \begin{bmatrix} 0 & 0^{(8nl_i+10)\times 2} & 0 \\ 0 & diag(a_i\ 0) & (0\ b_i)^T \\ 0 & (0\ b_i) & 0 \end{bmatrix}$ and $\mu_i = \begin{bmatrix} x_i \\ {}_Q x_i(x_i) \\ \mu_{i,end} \end{bmatrix}$

s.t.

1. $\mu_i^T \Pi_i^{\breve{S}_i} \mu_i + \breve{d}_i = 0;\ \mu_i^T \Pi_i^{\bar{S}_i} \mu_i + \bar{d}_i = 0;$ real and reactive power balance
2. $\mu_i^T \Pi_i^{\breve{F}_{i-j}^i} \mu_i = 0;\ \mu_i^T \Pi_{i-j}^{\bar{F}_{i-j}^i} \mu_i = 0;\quad j \in Br\{i\};$ definitions for the real and the reactive power flows as shown in 2 (3)
3. $\mu_i^T \Pi_{i-j}^{F_{flow}^{\max}} \mu_i - cap_l^2 \le 0;$ thermal flow limit
4. $\mu_i^T \left(\Pi_m^{\breve{G}_m} + \mathbf{j}\Pi_m^{\bar{G}_m}\right)\mu_i + p_{m,i}^{\max} p_{m,i}^{\min} + \mathbf{j}q_{m,i}^{\max} q_{m,i}^{\min} \le 0, m \in G\{i\};$ real and reactive power generation limits
5. $E_i^{\min} \le \mu_i^T \Pi_i^{E_i} \mu_i \le E_i^{\max};$ limits for the voltage magnitudes
6. $A_j^T \mu_i - \Phi_L^T v_L - \Phi_M^T v_M = 0 \leftrightarrow x_i - \Phi_L^T v_L - \Phi_M^T v_M = 0;$ nodal variables in terms of the central voltages where $A$ is a block matrix to collect an element $\xi_j = A_\xi^T \mu_j$;

$\mu_i^{nx\times 1} = \left[\alpha_i^{4\times 1}; \beta_i^{4\times 1}; \gamma_{i-l,i}^{4nl_i\times 1}; \delta_{i-l,i}^{4nl_i\times 1}; \omega_i^{2\times 1}; \breve{f}_{i-l,i}^{nl_i\times 1}; \bar{f}_{i-l,i}^{nl_i\times 1}; g_i^{2ng_i\times 1}; \mu_{i,end}^{1\times 1}\right]$; $\left(\alpha_i^T\ \beta_i^T\right)^T = A_{\alpha\beta}^T \mu_i$; $\left(g_i^{pT}\ g_i^{qT}\ \mu_{i,end}^T\right)^T = A_g^T \mu_i$ $\left(\breve{f}_{i-l,i}\ \bar{f}_{i-l,i}\right)^T = A_{fl}^T \mu_i$

; $\left(g_{m,i}^p\ g_{m,i}^q\ \mu_{i,end}\right)^T = A_{gm}^T \mu_i$ ; $\omega_i = A_\omega^T \mu_i$ ; $A_j = \begin{bmatrix} A_{\alpha\beta} & A_{\gamma\delta_1} & \cdots & A_{\gamma\delta_{Br\{j\}}} \end{bmatrix}$ ; $\Pi_{i-l,i}^{F_{flow}^{\max}} \triangleq A_{fl} A_{fl}^T$ ; $\Pi_i^{E_i} \triangleq A_\omega A_\omega^T$ ;

$\Phi_L = \left[\Phi_{\breve{S}_i} \Phi_{\bar{S}_i} \Phi_{\breve{F}_{i-1}^i} \Phi_{\bar{F}_{i-1}^i} \cdots \Phi_{\breve{F}_{i-Br\{i\}}^i} \Phi_{\bar{F}_{i-Br\{i\}}^i}\right]$ ; $\Phi_M = e_j + e_{j+Nb}$ ; $\Pi_j^{S_{inj}} \triangleq A_{\alpha\beta} BD \begin{pmatrix} \Pi_j^{S_j} \\ \mathbf{j}\Pi_j^{\bar{S}_j} \end{pmatrix} A_{\alpha\beta}^T - \frac{1}{2} A_g \begin{bmatrix} 0 & 0 & 1^{ng_j} \\ 0 & 0 & \mathbf{j}1 \\ 1^T & \mathbf{j}1^T & 0 \end{bmatrix} A_g^T$ ;

$\Pi_{i-l,i}^{F_{flow}} \triangleq A_{\gamma\delta_l} \begin{bmatrix} \Pi_{i-l,i}^{\breve{F}_{i-j}^i} & 0 \\ 0 & \mathbf{j}\Pi_{i-l,i}^{\bar{F}_{i-j}^i} \end{bmatrix} A_{\gamma\delta_l}^T - \frac{1}{2} A_{fl} \begin{bmatrix} 0 & 0 & 1 \\ 0 & 0 & \mathbf{j}1 \\ 1^T & \mathbf{j}1^T & 0 \end{bmatrix} A_{fl}^T$ ; and $\Pi_{m,i}^{G} \triangleq A_{gm} \begin{bmatrix} 1 & 0 & -p_{m,i}^{av} \\ 0 & \mathbf{j} & -\mathbf{j}q_{m,i}^{av} \\ -p_{m,i}^{av} & -\mathbf{j}q_{m,i}^{av} & 0 \end{bmatrix} A_{gm}^T$ ;

$p_{m,i}^{av} = \left(p_{m,i}^{\max} + p_{m,i}^{\min}\right)/2$, and $q_{m,i}^{av} = \left(q_{m,i}^{\max} + q_{m,i}^{\min}\right)/2$.

$$\min_{x_i \in X_i} \sum_i w_i(x_i, y, z_i) = \min_{x_i \in X_i} \sum_i \left[f_i(x_i) + g_i(x_i, y, z_i)\right] \tag{4B}$$

where $y = \left(v_L^T\ \ v_M^T\right)^T$ , $f_j(x_j) \triangleq \min_{x_i \in X_i} \sum_i x_i^T M_{ob}^i x_i$ , $g(x_j, y, z_j) = \frac{\rho_j}{2}\left\|x_j - \Phi_j^T y\right\|^2 + z_j^T\left(x_j - \Phi_j^T y\right)$ , or equivalently $g(x_j, y, z_j) = \frac{\rho_j}{2}\left\|x_j - \Phi_j^T y + z_j/\rho_j\right\|^2 - \frac{1}{2\rho_j}\left\|z_j\right\|^2$; and $X_j$ is the column space of $\Phi_j$. at Node $j$ defined by the constraints listed in 1–5 in (4A).

Note that all $A$'s are full column rank matrices that collect relevant parts from $\mu_i$, and that $\Phi_i$ includes the linear relationships between the local variable $x_i$ and the central variable $v$ and $y$. Using the definition of the function $f_j$, the following observations are made: 1) $f_j$ is a nonconvex, but smooth, function that is $C^1$ on an open set containing $X_j$ (defined by the column space of $\Phi_j$), and $\nabla f_j$ is Lipschitz continuous on $X_j$; 2) $X_j$ is nonempty, closed, and convex; and 3) $w_j$ is coercive.

Even though Eqs. (4A) and (4B) have low cardinalities in the decision variables due to the nonconvex nature of the problem, the uniqueness and existence of the solution are not guaranteed. To address the complexity issue, a surrogate function $h_j$ approximating Eq. (4B) is defined as follows:

$$H(x,y,z)=\sum_j h_j\left(x_j,y,z_j\right)=\sum_j\left[u_j\left(x_j\right)+g_j\left(x_j,y,z_j\right)\right] \text{ where } u_j\left(x_j\right)=\text{SDP of } f_j\left(x_j\right) \tag{5}$$

Because $u_j$ is the SDP relaxation of $f_j$, it is convex, Lipchitz continuous, and continuously differentiable on $X_j$. Note that $u_j$ relaxes the rank constraint only; therefore, the first derivatives of $u_j$ and of $f_j$ regarding $x_j$ are identical (i.e., $\nabla_{x_j} f_j\left(x_j\right)=\nabla_{x_j} u_j\left(x_j\right)$). Because the SDP is convex, the solution at the $k^{th}$ iteration is uniquely determined by $\tilde{x}_j^{k+1} \triangleq \arg\min_{x_j\in X_j} h_j\left(x_j,y^k,z_j^k\right)$ at given $y^k$ and $z_j^k$. Note that the relaxed nodal OPF is convex and has a small number of variables; therefore, its computational complexity remains manageable.

# Properties of $u_j$, $f_j$, $g_j$, $h_j$, $w_j$. and $\tilde{x}_j$

The nodal problem $f_j$ is $C^\infty$ smooth, but nonconvex, and prox-regular at $x_j$ relative to $f_j^\nabla$ if $x_j \in \text{dom}(f_j)$, $f_j^\nabla \in \partial f_j\left(x_j\right)$. Further, whenever $\left\|\bar{x}_j-x_j\right\|<\varepsilon_j^{prox}$ and $\left\|\hat{x}_j-x_j\right\|<\varepsilon_j^{prox}$, there exists $\varepsilon_j^{prox}, r_j^f>0$ :

$$f1)\, f_j\left(\bar{x}_j\right)\geq f_j\left(\hat{x}_j\right)+\left(f_j^\nabla\right)^T\left(\bar{x}_j-\hat{x}_j\right)-\frac{r_j^{\nabla f}}{2}\left\|\bar{x}_j-\hat{x}_j\right\|^2 \text{ ; and } f2)\, f_j\left(x_j^{k+1}\right)\leq f_j\left(x_j^k\right)+\nabla_{x_j} f_j\left(x_j^k\right)^T\left(x_j^{k+1}-x_j^k\right)+\frac{d_j^f}{2}\left\|x_j^{k+1}-x_j^k\right\|^2 .$$

Statement *(f2)* holds by the Descent Lemma [23] with the properties of $g_j$ and of $h_j$ described below. Its SDP relaxation $u_j$ is strongly convex, Lipschitz continuous, and prox-regular, and its first derivative equals that of $f_j$ (i.e.: $u1)\ u_j\left(\bar{x}_j\right)\geq u_j\left(x_j\right)+\left[\nabla_{x_j} u_j\left(x_j\right)\right]^T\left(\bar{x}_j-x_j\right)+\frac{c_j^u}{2}\left\|\bar{x}_j-x_j\right\|^2$ $u2)\ u_j\left(\bar{x}_j\right)\leq u_j\left(x_j\right)+L_j^u\left\|\bar{x}_j-x_j\right\|$ ; and $u3)\ \nabla_{x_j} u_j\left(x_j\right)=\nabla_{x_j} f_j\left(x_j\right)$. Whereas the variables considered in $f_j$ and $u_j$ are all local primal variables only, ADMM constraint $g_j$ depends on the central variable $y$ and the multipliers $z$. With the properties of $u_j$, $f_j$, and $g_j$, $w_j$ is $C^\infty$ smooth, but nonconvex, and prox-regular, whereas $h_j$ is strongly convex and Lipschitz continuous. Because all $X_j$ are bounded sets, $w_j$ and $h_j$ are coercive and lower-bounded over $X_j$. Because $x_j^{k+1}$ is the optimizer of convex $h_j$ (i.e., $\tilde{x}_j^{k+1} \triangleq \underset{x_j \in \tilde{X}_j;\, x_j \in X_j}{\arg\min}\, h_j\left(x_j, y^k, z_j^k\right)$), the following holds:

$$x1)\ \nabla_{x_j} u_j\left(x_j^{k+1}\right) + z_j^k + \rho_j\left(x_j - \Phi_j^T y\right) = 0 \text{ and } x2)\ h_j\left(x_j^{k+1}, y^k, z_j^k\right) - h_j\left(x_j, y^k, z_j^k\right) \le 0\ \forall x_j \in X_j.$$

# Proposed algorithm with provable convergence

## Overview of proposed algorithm

Even though the number of variables in the nodal OPF is much smaller than that in the original problem, the complexity of the problem makes it difficult to solve the nodal problem due to its nonconvex nature. To manage the nonconvex nature of the problem, the nonconvex components are relaxed. Once identified, the relaxed solution is mapped onto the feasible region, to the closest point in the feasible region from the identified solution. There are three cases: 1) the convexified problem is not feasible, 2) the convexified problem is feasible, but the solution to the problem is too far from the feasible region, and 3) the solution to the convexified problem is sufficiently close. Only if the mapped point is close enough, the nodal variables are updated and fed to the center node to update voltages.

## Proposed algorithm

We propose the following algorithm to solve $w_j$ in conjunction with $h_j$ in Eq. (5):

---

*Distributed, Regulated, Optimally-Homogeneous, and Scalable (DROHS) Algorithm*

---

1. Set $k = 0$, initialize all parameters, such as $\Delta^k \in (0,1)$, $y^k$, $\rho_j$, $\left\{z_j^k\right\} \in Null\left(\Phi_j\right)$.

2. If $x_j^k$ satisfies the termination criteria, terminate the algorithm.

3. Solve the local problems, $\tilde{x}_j^{k+1} \triangleq \underset{x_j \in X_j}{\arg\min}\, h_j\left(x_j, y^k, z_j^k\right)$.

4. Depending on the solution $x_j^{k+1}$, determine $\zeta_j^{k+1}$ (Fig 3)

- Case 1: $x_j^{k+1}$ is in the real space of convex $X_j$, i.e., in the column space of $\Phi_j$, $\zeta_j^{k+1} = x_j^{k+1}$.

- Otherwise $\left(x_j^{k+1} \notin X_j\right)$: Find a feasible projection of $x_j^{k+1}$ onto the feasible region $X_j$, $\zeta_j^{k+1}$near $x_j^{k+1}$;

- Case 2: the projection $\zeta_j^{k+1}$ is close enough to $x_j^{k+1}$ (i.e., $\left\|\zeta_j^{k+1} - x_j^{k+1}\right\| \le \varepsilon_j^k$ where $\varepsilon_j^k \triangleq \tau_j^k \left\|\tilde{x}_j^{k+1} - x_j^k\right\|$); accept the solution $\zeta_j^{k+1}$, $\zeta_j^{k+1} = \zeta_j^{k+1}$;

- Case 3: $\zeta_j^{k+1}$ is not sufficiently close; reject the solution $\zeta_j^{k+1}$, $\zeta_j^{k+1} = x_j^k$;

- update $\hat{x}_j^{k+1}$, $\hat{x}_j^{k+1} = x_j^k + \Delta^k\left(\zeta_j^{k+1} - x_j^k\right)$;

- set the error bound asymptotically to vanish as iteration proceeds (i.e., $\lim_{k\to\infty} \tau_j^k = 0$).

5. Compute $y^{k+1} = \arg\min_y h\left(\hat{x}_j^{k+1}, y, z_j^k\right)$, and update $z_j^{k+1} = z_j^k + \rho_j\left(\hat{x}_j^{k+1} - \Phi_j^T y^{k+1}\right)$ accordingly; and $\Delta^{k+1} = \Delta^k - a\left(\Delta^k\right)^2$ where $a \in (0.5, 1)$.

6. Update $x_j^{k+1}$ so that all $x$-variables are consistent with global variable $y$ (i.e., $x_j^{k+1} - \Phi_j^T y^{k+1} = 0$).

7. $k \leftarrow k + 1$, and go to 2.

---

Fig 3 illustrates how the $x$-optimization process determines the $\hat{x}_j^{k+1}$ described in Rule 4 of the algorithm. If the solution $x_j^{k+1}$ is feasible in $X_j$, $\hat{x}_j^{k+1}$ is the solution (Case 1). If the solution $x_j^{k+1}$ is not in $X_j$ but the point projected onto the convex region $X_j$ is sufficiently close, $\hat{x}_j^{k+1}$ is the projected point (Case 2). For Case 3, where either no solution is identified or the identified solution is not in $X_j$, and the distance to the projected point is too far away, the solution is rejected, and $\hat{x}_j^{k+1}$ is the point determined in the previous step. Upon the determination of $\hat{x}_j^{k+1}$, the iteration proceeds to a point between the current and the projected point depending on parameter $\Delta^k$.

After the $x$-optimization is performed, $y$-optimization follows. The $y$-optimization at given $\hat{x}_j^{k+1}$ and $z^k$ is a least square problem – convex problem. It should be noted that the solution to the convex problem is always feasible. Rule 6 is a striking feature of the proposed algorithm that is critically different from the ADMM approaches. Instead of the locally independent update of the primal variables in the ADMM approaches, the nodal variables are updated 1) linearly for the computational efficiency, and 2) with respect to the central variable so that the update is reasonably agreeable among nodal variables. The proof of the convergence of the proposed algorithm is presented in this paper's Appendix B.

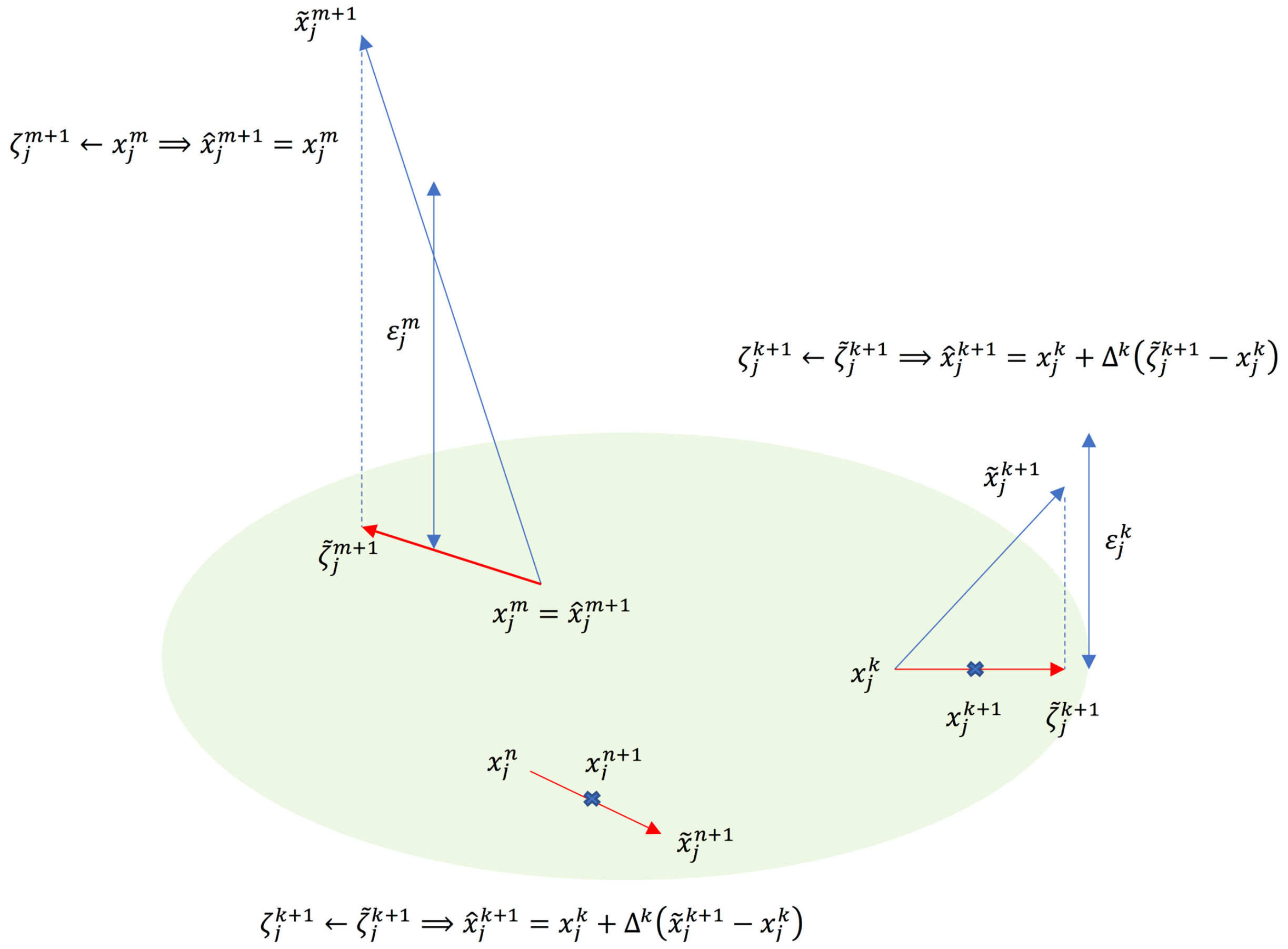

Fig 3. $x$-update from the result of optimization with 1) a feasible solution $\tilde{x}_j^{n+1}$, 2) an infeasible solution, but the distance between the solution and the projection onto the feasible region is close enough $\left\|\tilde{\zeta}_j^{k+1} - \tilde{x}_j^{k+1}\right\| \le \varepsilon_j^k$, and 3) an infeasible solution in that the projection onto the feasible region is too far $\left\|\tilde{\zeta}_j^{m+1} - \tilde{x}_j^{m+1}\right\| \le \varepsilon_j^m$.

# Implementation Details

## Implementation of Rule 4

The local SDP relaxation yields a matrix (known as $Z_j^{k+1}$) instead of a vector $x_j^{k+1}$; therefore, a new criterion should be established to compute $\zeta_j^{k+1}$, as follows: Because $Z_j^{k+1}$ is from SDP, all the eigenvalues are non-negative (i.e., $Z_j^{k+1} = \zeta_j^{k+1}\left(\zeta_j^{k+1}\right)^T + \sum_{m=2}^{r} \psi_{j,m}^{k+1}\left(\psi_{j,m}^{k+1}\right)^T$) where $\psi_{j,m}^{k+1}$ is the $m^{\text{th}}$ column vector, and $r$ is the rank of $Z_j^{k+1}$:

$$\left\|\zeta_j^{k+1}\left(\zeta_j^{k+1}\right)^T - Z_j^{k+1}\right\| \le \left\|\zeta_j^{k+1} + x_j^{k+1}\right\|\left\|\zeta_j^{k+1} - x_j^{k+1}\right\| \tag{6}$$

To construct an equivalent criterion to $\left\|\zeta_j^{k+1} - x_j^{k+1}\right\| \le \varepsilon_j^k$ with $\left\|\zeta_j^{k+1} + x_j^{k+1}\right\| \le 2\left\|\zeta_j^{k+1}\right\| + \left\|\zeta_j^{k+1} - x_j^{k+1}\right\| \le 2\left\|\zeta_j^{k+1}\right\| + \varepsilon_j^k$, we impose $\left\|\zeta_j^{k+1}\left(\zeta_j^{k+1}\right)^T - Z_j^k\right\| \le \varepsilon_j^k\left(\varepsilon_j^k + 2\left\|\zeta_j^{k+1}\right\|\right)$ to determine whether or not the feasible projection is acceptable. The grouping criterion is simplified to $\frac{1}{2\left\|\zeta_j^{k+1}\right\|}\left\|\sum_{m=2}^{r}\psi_{j,m}^{k+1}\left(\psi_{j,m}^{k+1}\right)^T\right\| \le \varepsilon_j^k$ by invoking $\varepsilon_j^k \ll \left\|\tilde{\zeta}_j^{k+1}\right\|$. Let $\lambda_{Z_j^k}^m$ be the $m^{\text{th}}$ eigenvalue of $Z_j^k$. Then, the criterion becomes $\lambda_{Z_j^k}^2 \le 2\lambda_{Z_j^k}^1 \varepsilon_j^k$.

The definition of $\varepsilon_j^k \triangleq \tau_j^k\left\|\tilde{x}_j^{k+1} - x_j^k\right\|$ is also modified as $\varepsilon_j^k \triangleq \tau_j^k\left[\sqrt{\left\|x_j^k\right\|^2 + \left\|\tilde{Z}_j^{k+1} - x_j^k\left(x_j^k\right)^T\right\|} - \left\|x_j^k\right\|\right]$ using $\left\|Z_j^{k+1} - x_j^k\left(x_j^k\right)^T\right\| \le \left(2\left\|x_j^k\right\| + \left\|x_j^{k+1} - x_j^k\right\|\right)\left\|x_j^{k+1} - x_j^k\right\|$.

# Choice of parameters

In contrast to the ADMM approaches, only a few parameters are assigned in the proposed algorithm. Because $\left\{z_j^k\right\} \in Null\left(\Phi_j\right)$, and the specific choice of $z_j$ is irrelevant, all $z_j$'s are random vectors in the null space of known $\Phi_j$. It is important to hold $\sum_{k=1}^{\infty}\left\|\varepsilon^k\right\|^2 < \infty$ to guarantee the convergence of the proposed algorithm [24]. Due to the fact that $\sum_{k=1}^{\infty}\left\|\varepsilon^k\right\|^2 \le \left\|x_j^{k+1} - x_j^k\right\|_{\max}^2 \sum_{k=1}^{\infty}\left(\tau_{\max}^k\right)^2$, the maximum $\tau$ over all nodes at the $k^{\text{th}}$ iteration is set by $\tau_{\max}^k = \tau_{\max}^0 / k$. The parameters used were as follows: $\varrho_j = 20$ for $v_L$ and 200 for $v_M$; $\Delta^0 = 0.3$; and the number of maximum iteration = 100; $\tau_{\max}^0 = 10^{-3}$.

The global variable $y$ is randomly assigned (cold start) to test the algorithm's robustness. For comparison, a flat start and a warm start are also attempted. Here, a flat start refers a point at which real voltages are unity, and imaginary voltages are within [-0.1, 0.1]; and a warm start [20] is a solution to the PF. It is recognized that a faster convergence is obtained when the initial point is close to the feasible region, which is commonly observed in numerical iterative methods.

## Termination criterion

The algorithm is terminated when no further progress is made. The progress is measured in terms of the global variable $y$, and the local variables $x$ and $z$. After the solution is identified, the objective function $W^*$ is determined according to (4B) where $W(x,y,z) \triangleq \sum_i w_i(x_i,y,z_i)$ , i.e., $W^* = W(x^*,y^*,z^*) = \sum_i w_i^*(x_i^*,y^*,z_i^*)$. The interim value of W at the $k^{th}$ iteration is $W^k = W^k(x^k,y^k,z^k) = \sum_i w_i^k(x_i^k,y^k,z_i^k)$, and the progress measure $\eta^k$ at the $k^{th}$ iteration is defined:

$$\text{Progress measure} \quad \eta^k \triangleq \frac{W^k - W^*}{W^*} \tag{7}$$

The criterion used in Rule 2 terminates the algorithm if the progress measure is less than the tolerance ($10^{-7}$) (i.e., $\eta^k \leq 10^{-7}$).

# Results and Discussion

## Simulation environments

The model systems used in the simulations are available from MATPOWER [25] for 3-, 4-, 9-, 14-, 24-, 30-, 39-, 57-, 85-, 118-, 300-, and 2,000-bus systems. All simulations were performed using a Mac pro with 2 × 2.93 GHz 6-core Intel Xeon processors and 6 GB 1333 MHz DDR3 memory. The local SDP problems were solved using the CVX solver [26]. To compare the results from various systems, the lines in the figures are normalized so that all start at zero. The initial points for all the cases are cold starting points, and all the real and imaginary voltages are set to random numbers (not even a flat start). The qualities of the solutions are all numerically identical to those identified using MATPOWER for the model systems tested (all the solutions $v^*$ are less than $10^{-4}$ from the MATPOWER solutions ($v_{MATPOWER}$), $\|v^* - v_{MATPOWER}\|/\|v^*\| \leq 10^{-4}$). For comparison, we also attempted a flat start and observed that $N_{iter}$ decreases at least 30%, but the flat start also finds the same solutions.

# Subsystems of the test cases

In the ADMM approaches referred to in the literature [14–20], the grid partitioning is performed, but the details of the partitioning are not well listed. In considering the computational coupling, the partitioning is based on spectral clustering [20,27] while each subsystem contains at least one generator [14,20]. The inclusion of a generator in a subsystem seems a natural choice to fulfill the power balance equality constraints (1 in (3)) within each substation. However, the inclusion of a generator does not serve the purpose well, because the generator may not be dispatched if its generation costs are excessively high. The spectral clustering is performed using two different algorithms: 1) unnormalized [27] and 2) normalized [28]. The algorithms yield different grid clustering. For example, the unnormalized algorithm results in two clusters for the IEEE 3-bus systems, while the normalized algorithm finds a single cluster for the same system [28]. Fig 4 illustrates the path length, maximum number of nodes in a subsystem, and the number of subsystems.

The lines are best-fit lines, showing the positive correlations of path length (red dotted line path length $\propto Nb^{0.30}$), the maximum number of nodes in a subsystem (green solid line $\propto Nb^{0.43}$), and the number of subsystems (blue broken line $\propto Nb^{0.76}$) based on the unnormalized algorithm (left plot). The results using the algorithm with a different normalization [28] yield path length $\propto Nb^{0.34}$, the maximum number of nodes in a subsystem $\propto Nb^{0.45}$, and the number of subsystems $\propto Nb^{0.66}$ (right plot). Even though the results from the two plots are not exactly the same, the positive correlations with the system size are clear. They affect the computational efficiency of ADMM [14,20]; path length affects the communication costs and $N_{iter}$ because the decision at each subsystem should be delivered to the rest of the system; the number of subsystems affects the number of computation cores and the communication costs unless PL equals one; and the maximum number of nodes affects the computation costs at each subproblem. Therefore, the decrease in the computation costs at each subsystem is achieved in exchange for the increased costs of communication among subsystems. Because the path length and number of subsystems increase rapidly with system size, the improvement in computational efficiency of the ADMM approach is questionable.

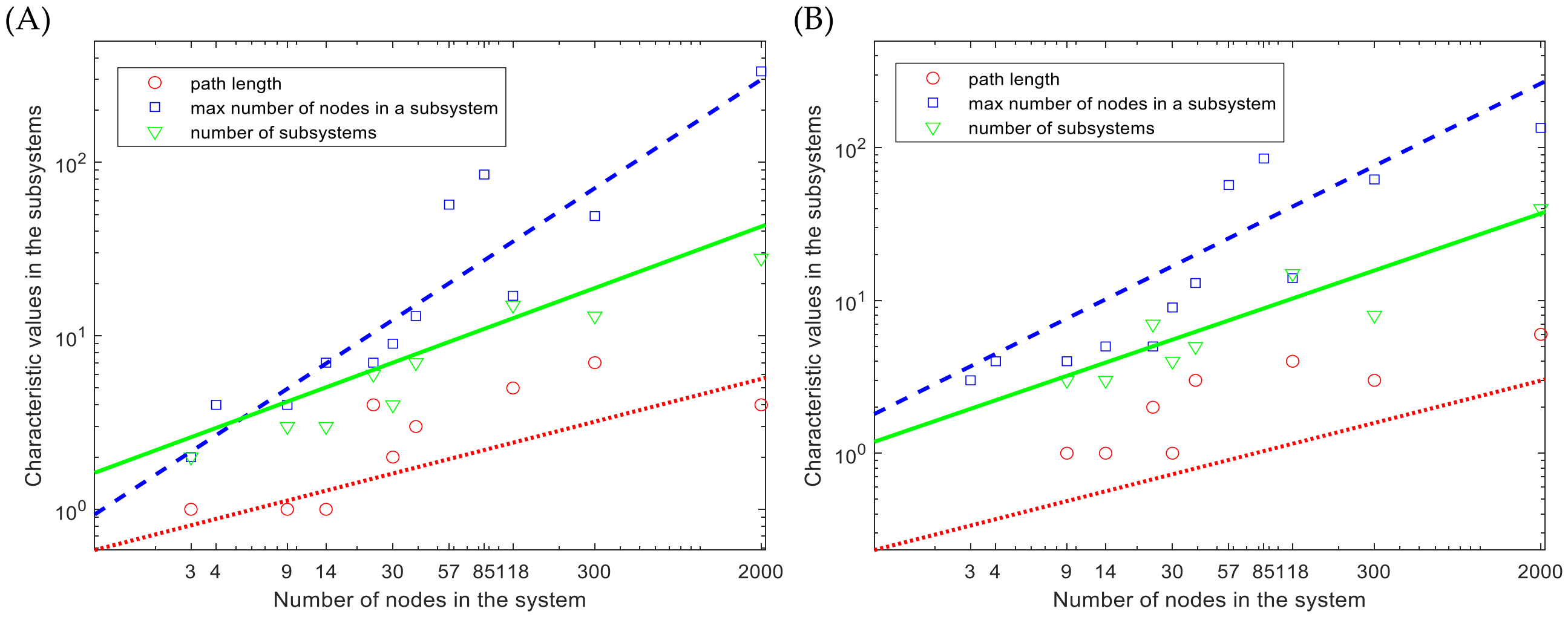

Fig 4. The path length, the maximum number of nodes in a subsystem, and the number of subsystems for the IEEE 3-, 4-, 9-, 14-, 24-, 30-, 39-, 57-, 85-, 118-, 300-, and 2,000-bus systems. The lines indicate the positive correlations on the size of the system. The red-dotted line, the blue broken line, and the green solid line represent the best-fit curves for the path length, the maximum number of nodes in a subsystem, and the number of subsystems, respectively. Two algorithms are applied for clustering nodes: (A) unnormalized algorithm, and (B) normalized based on the algorithm in [28].

# Maximum cardinality of nodal variables in the proposed algorithm

The nodal variables in the $x$-optimization $\mu_j$ are $\left[\alpha_i^{4\times1};\beta_i^{4\times1};\gamma_{i-l,i}^{4nl_i\times1};\delta_{i-l,i}^{4nl_i\times1};\omega_i^{2\times1};f_{i-l,i}^{nl_i\times1};\bar{f}_{i-l,i}^{nl_i\times1};g_i^{2ng_i\times1};\mu_{i,end}^{1\times1}\right]$. Because $\mu_{end}$ is unity, the cardinality of $\mu_j$ is 10$nl$ + 2$ng$ + 10 where $nl$ is the number of lines connected to the $j^{th}$ node, and $ng$ is the number of generators located at the node. In determining the voltages through the power channel, the node with a high value for $nl$ plays a key role. Fig 5 illustrates the maximum cardinality $n_{\text{var}}^{\max}$ of $\mu_j$ with the systems' sizes. The dashed line indicates $n_{\text{var}}^{\max} \propto Nb^{1/4}$. Even though the positive relationship between $n_{\text{var}}^{\max}$ and $Nb$ is visible among the model systems, the relationship is not necessarily positive. The number of variables in the central OPF is typically 2$Nb$ + 2$Ng$ where $Nb$ and $Ng$ are the numbers of buses and generators in the system, respectively. For small systems, such as 3-, 4-, 9-, 14-, and 24-bus systems, the maximum cardinalities of the nodal variables are higher than the number of variables in the central OPF; therefore, it would be challenging to keep the computational costs of the distributed OPF for the small systems lower than those of the central OPF. It is noteworthy that the communication costs remain manageable, because each node directly communicates with the central node due to the fact that PL equals unity.

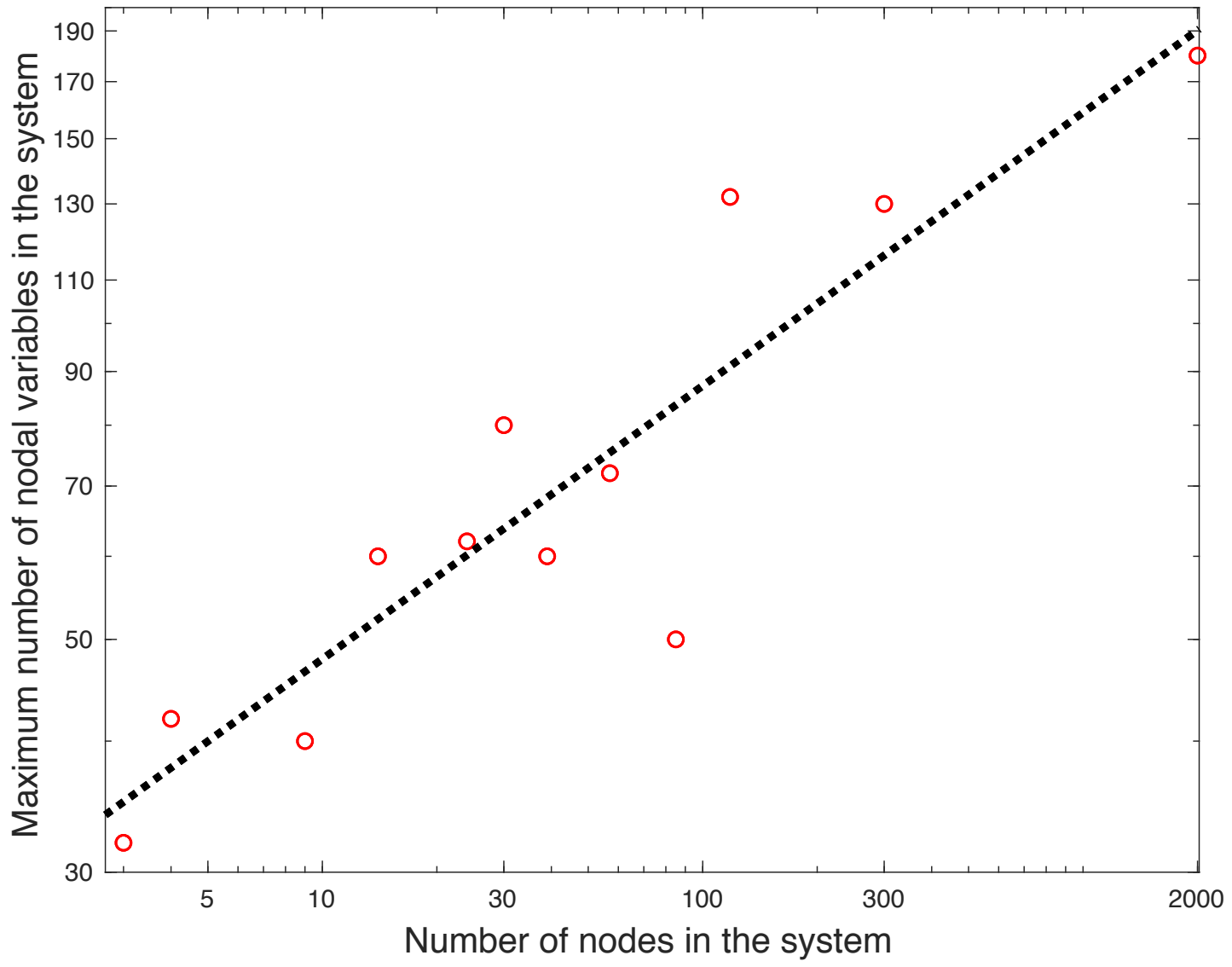

Fig 5. Cardinalities of nodal variables in $x$-optimization in terms of system size.

# Distributed OPF for 3-, 4-, 24-, 39-, and 85-bus systems

To examine whether the choice of parameters affects the convergences, optimizations were executed with many randomly assigned initial points, and uniform convergences were consistently observed (Fig 6). In general, $N_{iter}$ depends on the initial guess of the voltages. As the distance of the initial point from the solution becomes closer, $N_{iter}$ decreases, which is commonly observed in numerical iterative methods. All the solutions identified are numerically the same as the solutions using MATPOWER.

A worst-case convergence was observed for the IEEE 39-bus system. The nodal OPF for the system involves large negative eigenvalues (i.e., the nodal OPFs are highly nonconvex). According to Rule 4, a relatively large number of solutions are rejected because they are not close to the feasible regions and, therefore, $N_{iter}$ becomes large (approximately 40 iterations). However, the solution identified is a local minimizer, and $N_{iter}$ is still much smaller than the number of iterations for the ADMM approaches. For the visual presentations of the results obtained for various systems, the progress is normalized to the 1st iteration (i.e., the curves begin at 0 for the 1st iteration).

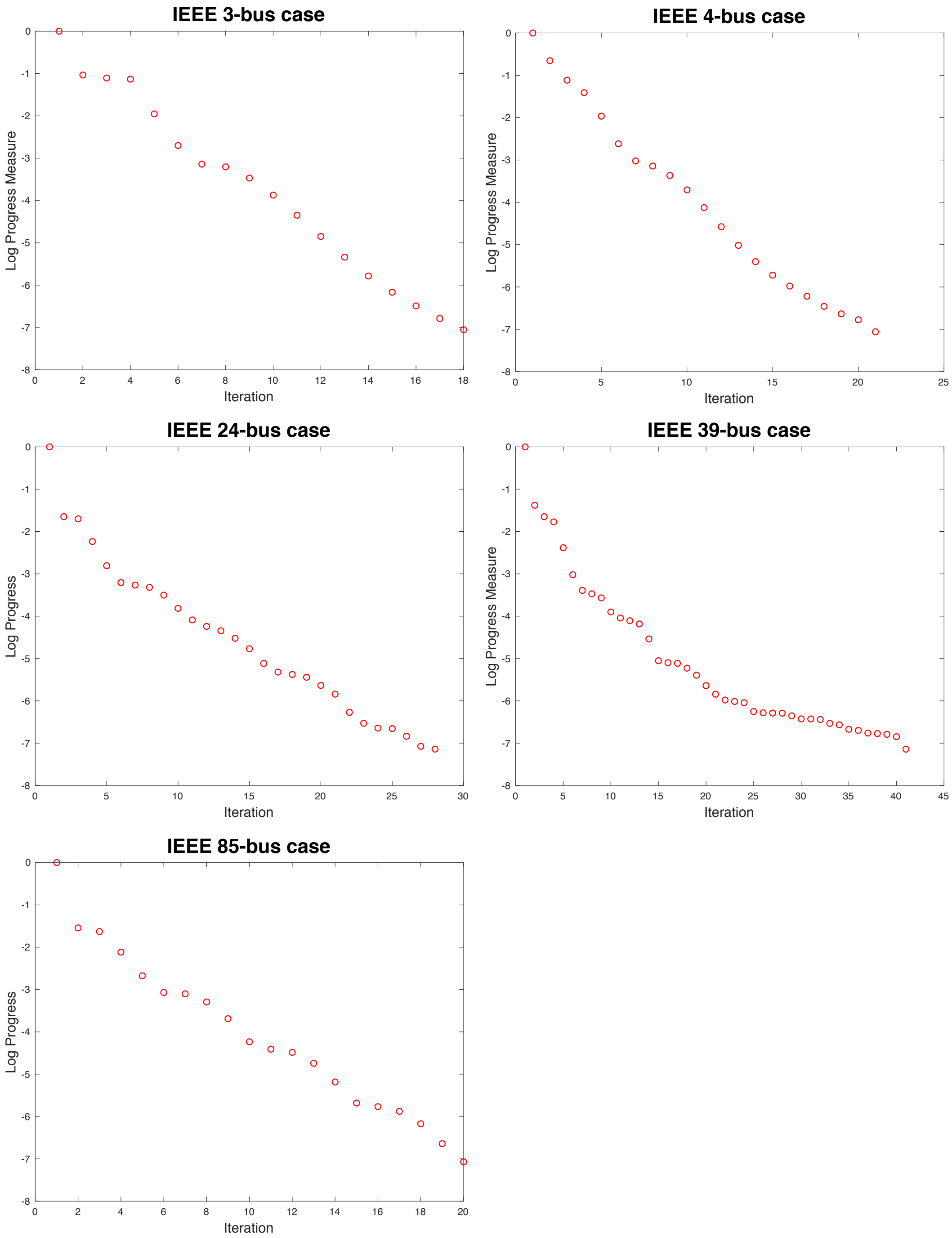

Fig 6. Convergence of proposed algorithm for selected systems (i.e., IEEE 3-, 4-, 9-, 14-, 24-, 39-, and 85-bus systems).

# Comparison of convergence to ADMM approaches

The convergence measures were compared to those from the ADMM approaches reported in multiple studies [14], [17], [18], and [19]. It is worth noting that Erseghe [14] and Zhang et al. [17] do not take the flow limits into consideration; the approach of Engelmann et al. [18] involves a high communication cost; and the approach of Madani, Kalbat, and Lavaei [19] may yield a physically infeasible solution. We attempted the ADMM approaches with the flow constraints and feasibility, but our implemented solver failed to converge with any starting points. Instead, we compared the convergence behaviors of the proposed algorithm to those reported in the previous [14,17–19]. For the visual presentation, the convergence behaviors are normalized so that all the curves begin at the same point (see Fig 7 for the comparison of the convergence).

In the previous research [14,17–19], the size of the subsystems increased with that of the system. Therefore, the tradeoff between communication costs amongst the subsystems and the computation costs for each subsystem make it difficult to develop a scalable algorithm. As the system size increases, $N_{iter}$ increases significantly for the results in two of the previous studies [14,19]. The performance measure in another study [17] fluctuates consistently with various systems, indicating that the convergence may not be guaranteed. Because the final study [18] requires the information exchange regarding sensitivities as well as the primal variables, the computation and communication costs per iteration increase much more rapidly than do the computation costs of a central OPF solver.

Different from the clustering containing multiple nodes used in the ADMM approach, the proposed algorithm for each subsystem contains only 1 node regardless of the system size; each subsystem directly communicates with the rest of the entire system; and the communication involves only the primal variables. With this difference in mind, the fast convergences observed in the proposed algorithm are similar to those in two of the previous studies [17,18], which are consistently much faster than those in the other two previous studies [14,19]. Whereas $N_{iter}$ are like those reported by Engelmann et al. [18], the convergences of the proposed algorithm occur more uniformly and consistently.

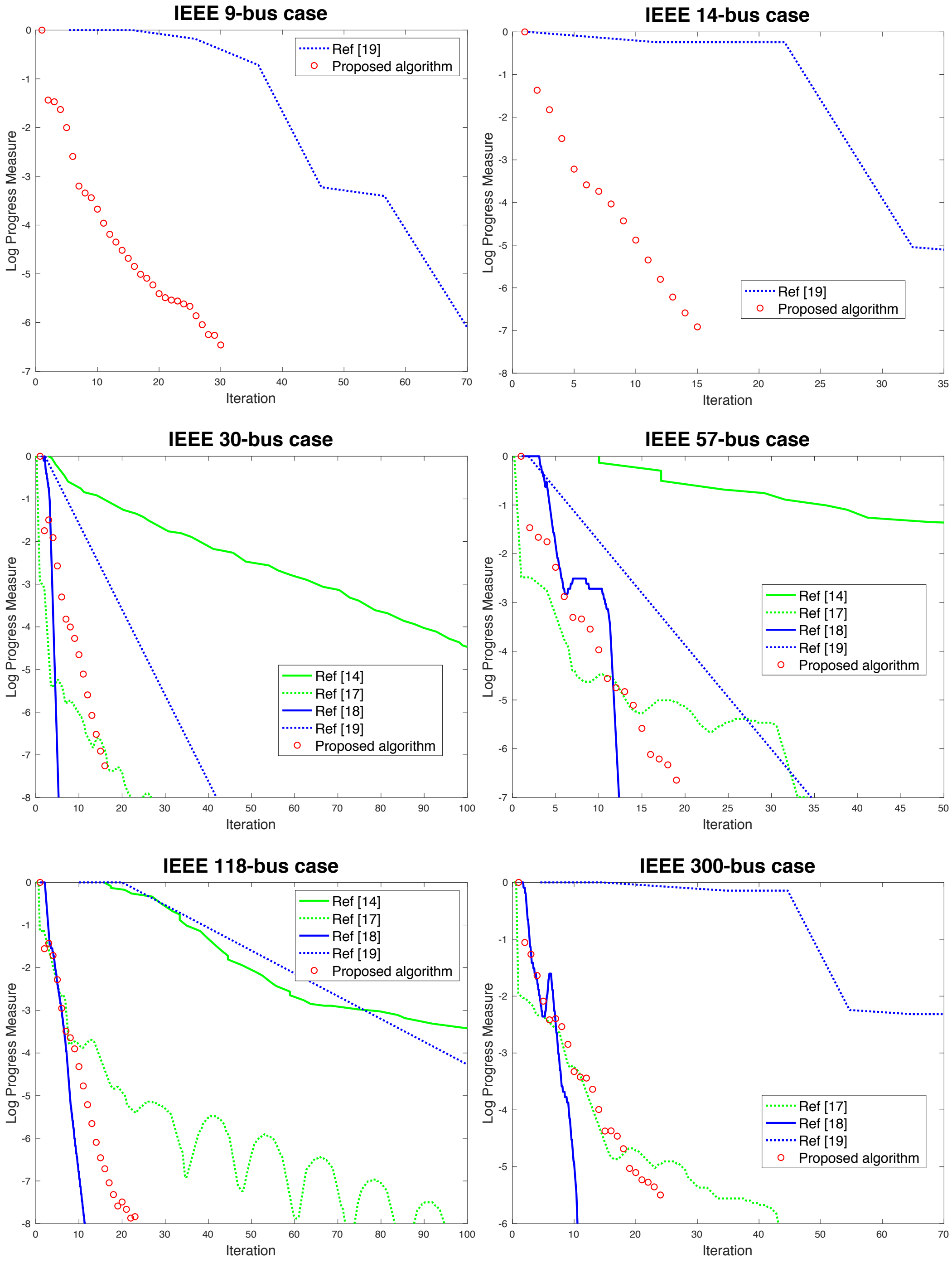

Fig 7. Convergence of proposed algorithm for selected small systems (i.e., IEEE 9-, 14-, 30-, 57-, 118-, and 300-bus systems). The convergences of the ADMM approaches [14, 17-19], are compared.

Many approaches using ADMM do not consider the flow limits [14,17], and/or feasibility [19] and, therefore, it is not appropriate to compare the quality of the solutions. Instead, we compared our solutions with those using a central OPF solver, MATPOWER. For all the cases described above, they yield numerically identical solutions. It is not clear why the proposed algorithm finds the same solution as the central OPF solver. We tested with the central SDP-relaxation for a small-scale system to estimate the "global" solution. Due to the memory issue, our SDP tests are limited up to 118-bus systems. When the relaxation returns the rank-1 solution, the solution is global and physically feasible. For several cases, the test cases are of the global solutions that are also identified by MATPOWER and by the proposed algorithm. However, there are cases where the SDP relaxation returns a physically infeasible solution. For these cases, both MATPOWER and the proposed algorithm identify numerically the same local minimizers that may not be global solutions. MATPOWER finds a point that meets the first-order necessary conditions for optimality [25]. Although it finds a minimizer in most cases, there is no guarantee that the identified solution is a minimizer—it can be either a maximizer or a saddle point. On the other hand, the proposed algorithm identifies a solution that meets the second-order optimality conditions, which guarantees that the solution is a minimizer [29]. From a practical point of view, there are two advantages of the proposed algorithms over MATPOWER: 1) they are numerically stable and, therefore, robust because all the subproblems are convex—no issues regarding the rank-deficient Hessian matrix, and 2) they use distributed computation and, therefore, manageable computation in each subproblem. However, there are disadvantages of the proposed algorithm: 1) the requirement of multiple cores to perform the distributed optimization, and 2) the number of nodal variables can be larger than that of the central OPF; for example, 3-, 4-, 9-, 14-, and 24-bus (See Fig 5).

## Large-scale OPF: 2,000-bus system

We tested the algorithm for a large-scale system, the 2,000-bus system that is a synthetic grid on a footprint of Texas. Fig 8 presents a convergence pattern. The solution identified is numerically identical to the one found using MATPOWER. $N_{iter}$ remains small, as is the case for the small systems (See Figs 6 and 7). Note that the nodal OPF includes only 1 bus, and the communication costs remain small because PL equals unity. If a sufficient number of computation cores are provided, the proposed algorithm is scalable if the computation cost for solving a subproblem does not increase significantly as the system size increases.

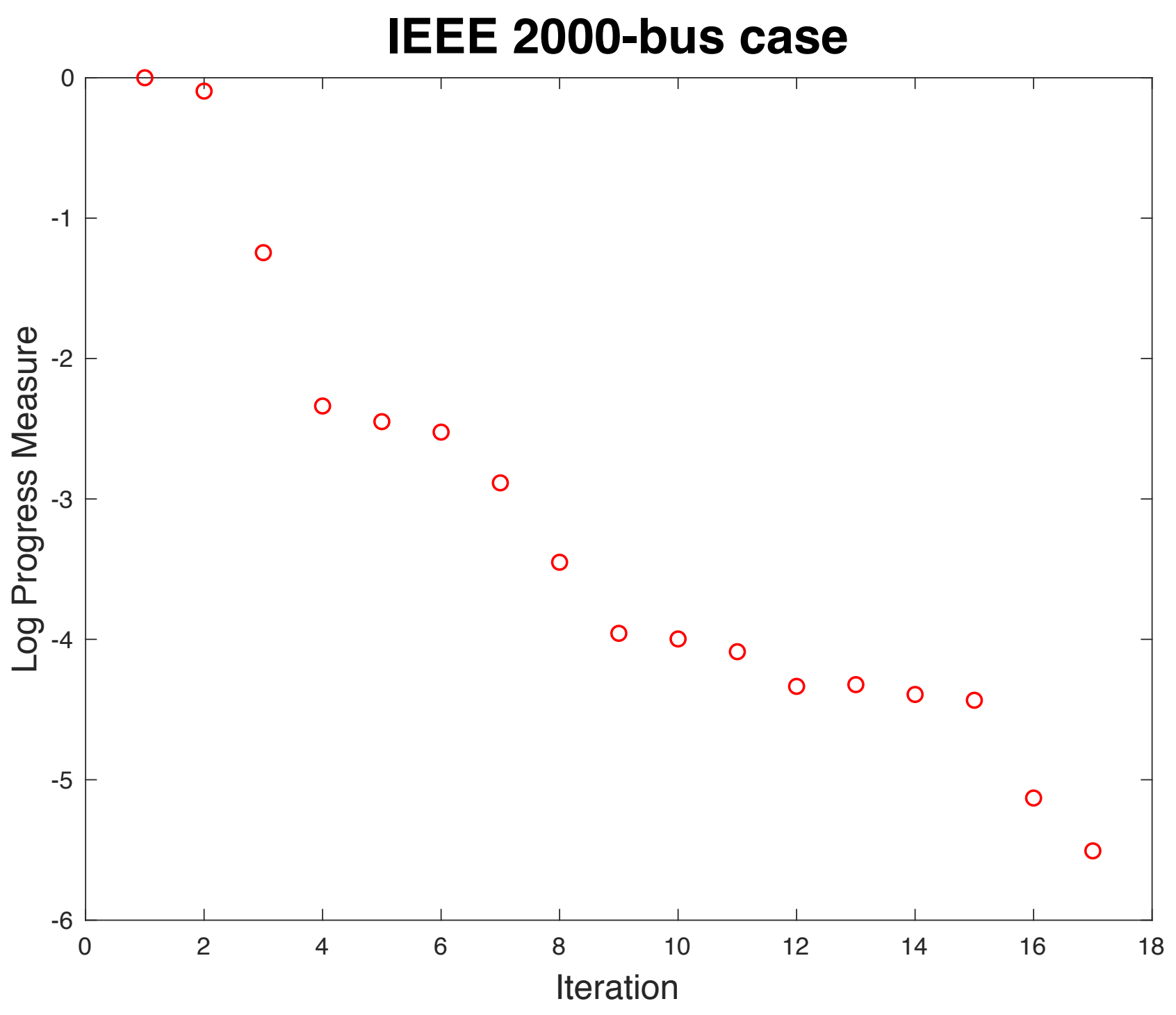

Fig 8. Convergence of proposed algorithm for synthetic 2,000-bus system on a footprint of Texas.

Different from other algorithms reported in the literature (i.e., [14–18]), the proposed algorithm converges regardless of the initial points, but the number of iterations decreases by at least a factor of 2 if it uses a flat start. Another shortcoming of the algorithms reported in previous research [14–18], is the quality of their solutions. Whereas the existing OPFs find low-voltage [17] or suboptimal [16] solutions, the proposed algorithm finds the same solution as a central OPF solver. One previous study [16] claimed that there might be an optimal number of partitioned subsystems due to the tradeoff between communications and computational costs used to obtain a local solution. The proposed method keeps the PL = 1, and the central computation is a simple addition of $x_j$ through both the power and voltage channels. From the simulations with various parameter values such as $\varrho_j$, $\Delta^k$, $\tau_j^1$, and $a$, we obtained the same solutions, indicating that the proposed algorithm is robust as well as efficient. In addition, the proposed model does not require any partitioning.

# Comparison to a heuristic central OPF solver, MATPOWER

In comparing the convergence between the proposed algorithm and the central OPF, the number of variables and $N_{iter}$ were examined. The maximum cardinality $n_{\text{var}}^{\max}$ of $\mu_{\text{j}}$ increases with the system size in $n_{\text{var}}^{\max} \propto Nb^{1/4}$ (See Fig 4). Fig 9 illustrates the number of variables in the central OPF as well as the maximum cardinality of the nodal variables. The blue line is the best-fit line for the central OPF ( $n_{\text{var}}^{OPF} = 2Nb + 2Ng \propto Nb^{0.96}$), whereas the red line is the best-fit line for the distributed algorithm. It is clear that $n_{\text{var}}^{OPF}$ increases with $Nb$ in a much faster way than $n_{\text{var}}^{\max}$ does. The black line is the boundary at which $n_{\text{var}}^{OPF}$ equals $n_{\text{var}}^{\max}$. If enough computation cores are available, the proposed algorithm involves reduced computation costs per iteration for systems larger than or equal to IEEE 24-bus. Whereas $n_{\text{var}}^{OPF}$ almost linearly increases with $Nb$ because $Ng \ll Nb$ in most systems, $n_{\text{var}}^{\max}$ does not necessarily increase theoretically. A key observation is that the proposed approach can be a scalable algorithm if $N_{iter}$ does not rapidly increase with the system size.

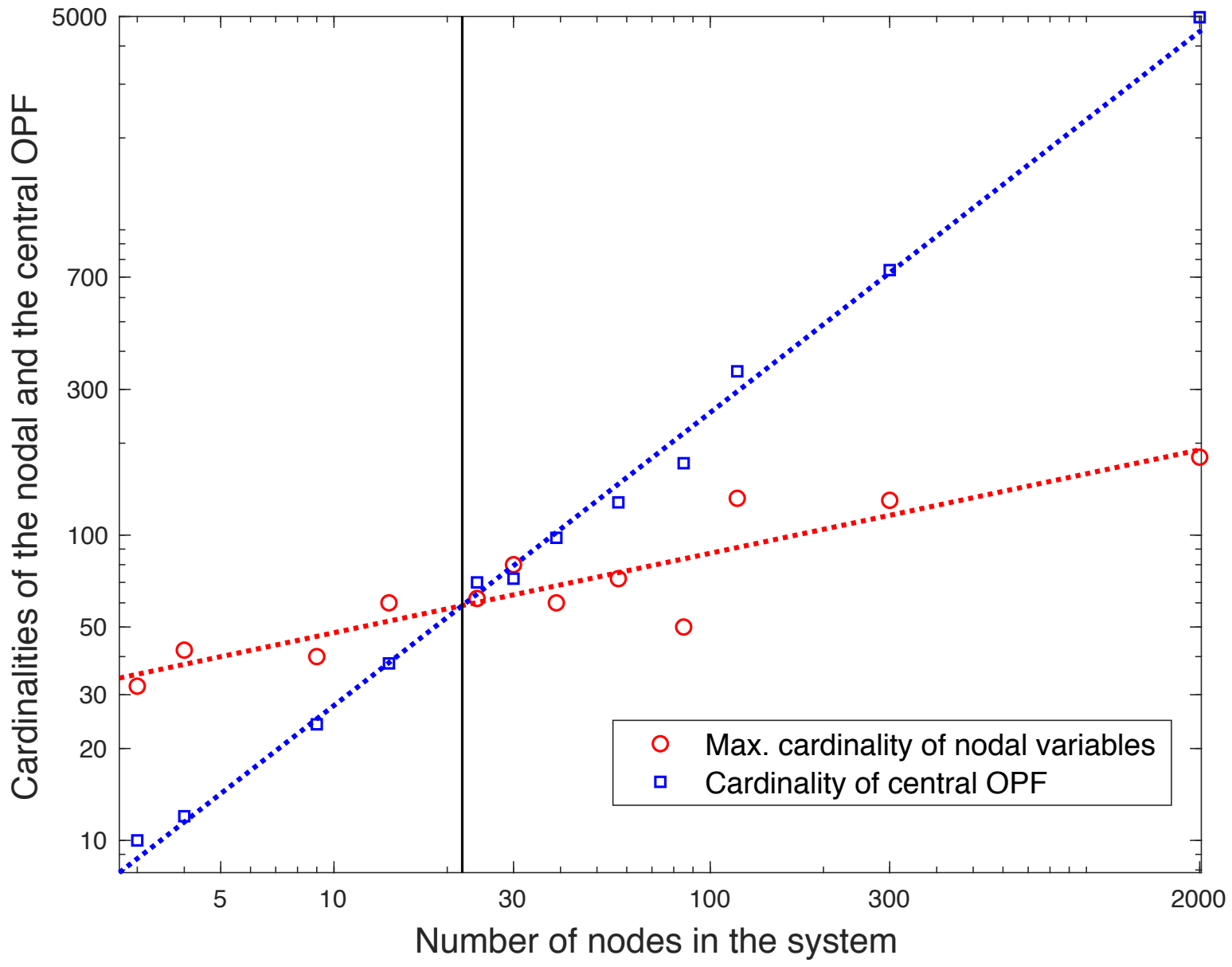

Fig 9. Cardinalities of central OPF and proposed algorithm for the test systems.

Fig 10 presents $N_{iter}$ for the central OPF ($Nit_{cent}$) and of the proposed algorithm ($Nit_{dist}$) on the tested systems. A visible increase in $Nit_{cent}$ is observed with $Nb$, but the dependence of $Nit_{dist}$ on $Nb$ is not evident—the solid lines are best-fit curves (log-log plot) that indicate $Nit_{cent} = Nb^{0.16}$ and $Nit_{dist} = Nb^{-0.007}$, respectively. The dotted lines are average $N_{iter}$, and the values are 15 and 23, respectively.

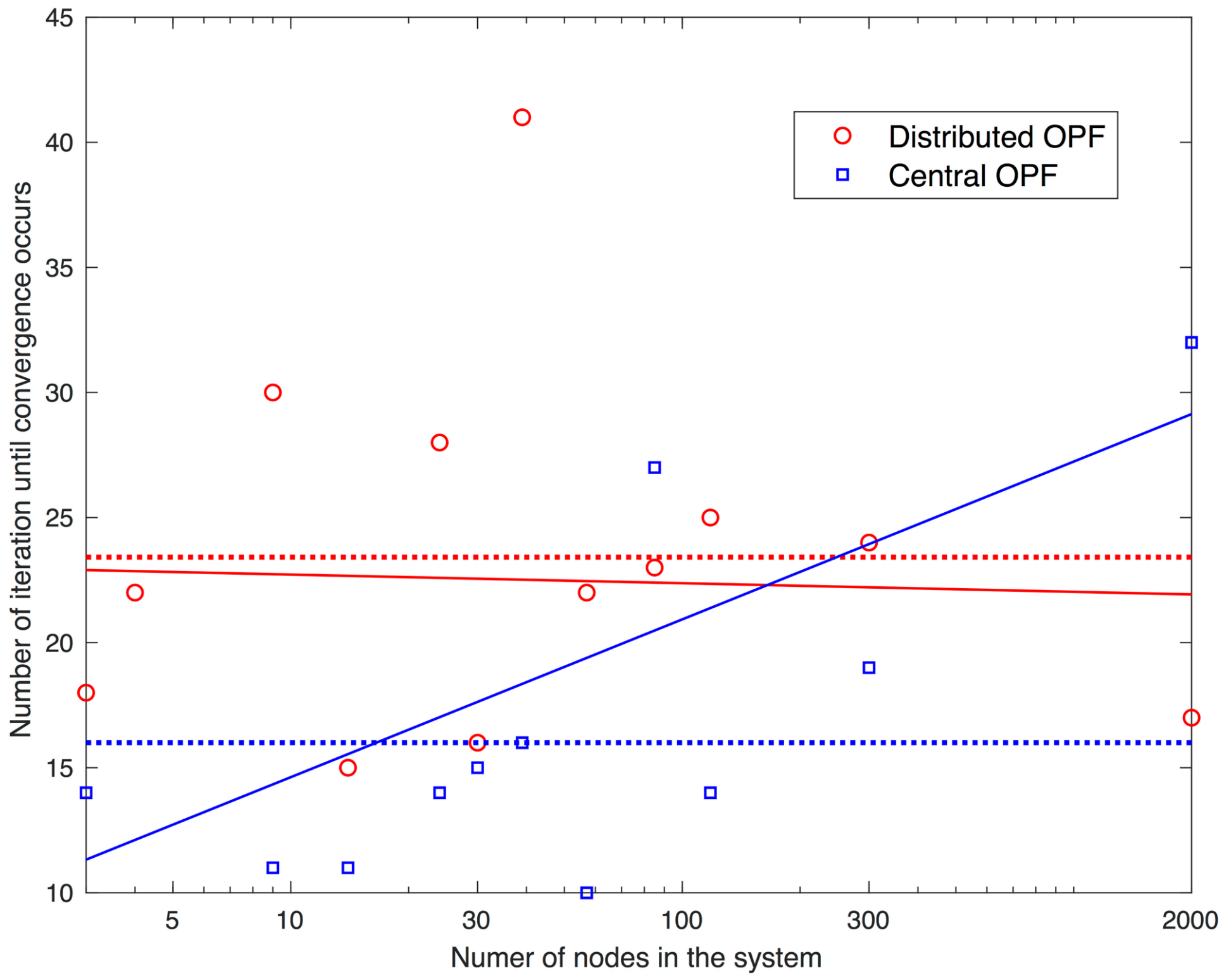

Fig 10. $N_{iter}$ for central OPF and of the proposed algorithm on the tested systems. Dotted lines are averages, and solid lines are best-fit curves.

From the comparisons (the number of variables and $N_{iter}$), we conclude that the computation cost of the distributed algorithm increases with $Nb$ at a much slower rate than that of the central OPF if sufficient computational resources are available.

# Scalability of the algorithm

The computation efficiency of a distributed computation depends on the number of iterations and the computational cost per iteration. The cost is determined by the local computation with the largest number of variables. Fig 10 shows that the number of iterations does not increase with system size.

Fig 11 illustrates the maximum nodal computation time depending on the maximum cardinality of the nodal variables. This performance dependence on the number of nodal variables may shed light on the claim by Loukarakis, Bialek, and Dent [30] that a larger system does not necessarily imply an inferior convergence performance. The dashed line indicates that the maximum nodal computation time is proportional to the maximum cardinality of the nodal variables with the power of 2.64. The cardinality to the power of 2.64 observed in this study is close to the theoretically estimated 3 for the SDP solver. Note that it is not necessary for the maximum cardinality of the nodal variables to be positively correlated with the system size. Rather, the maximum cardinality depends on the local grid topology of the system. For the tested systems, $n_{\mathrm{var}}^{\max} \propto Nb^{1/4}$ and $CT \propto \left(n_{\mathrm{var}}^{\max}\right)^{2.64}$, which results in $CT \propto Nb^{2.64/4} = Nb^{0.66}$ where $CT$ represents the maximum nodal computation time. The total computation cost is bounded by the product between $N_{iter}$ and CT. The total computation cost is bounded by $CT \propto \vartheta\left(\lceil Nb/N_{core} \rceil Nb^{0.66}\right)$ where $N_{core}$ is the number of cores. If the computational resource is sufficient, $CT \propto \vartheta\left(Nb^{0.66}\right)$. From this observation, the computation cost increases sub-linearly for a large-scale network.

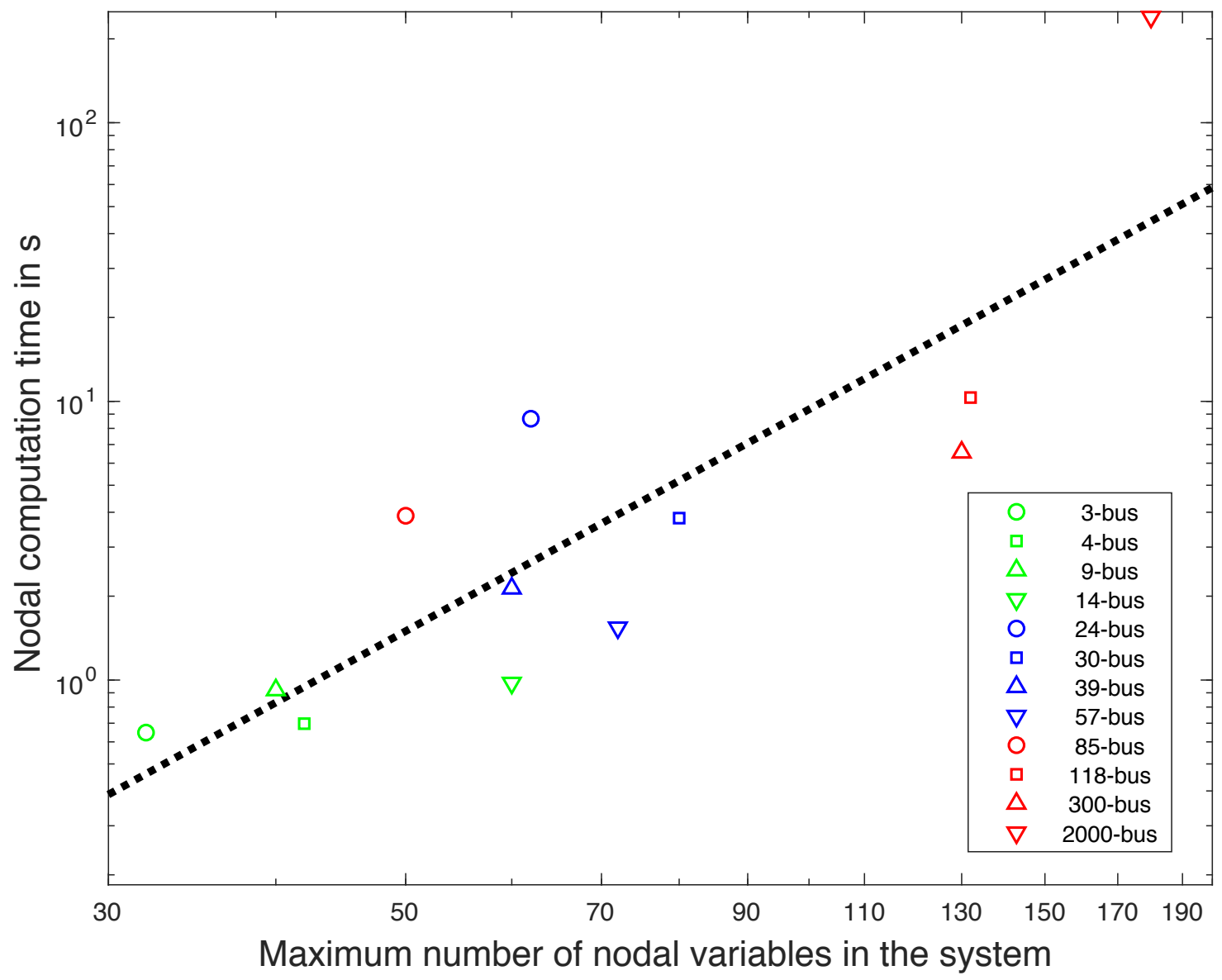

Fig 11. Computation times in seconds for nodal OPF of the maximum cardinality of the nodal variables. Dotted line is the best-fit line with the slope of 2.64.

If a single core is available for computation, the computation cost is in $\vartheta\left(Nb^{1+p_{SDP}/4}\right)$ where $p_{SDP}$ is the computational complexity for a nodal SDP. For the state-of-art central heuristic OPF solvers, the computation cost is in $\vartheta\left(Nb^{1.5}\right)$. Therefore, when the computation resources are highly limited, the proposed algorithm would still be efficient with a convex problem solver that yields $p_{SDP} \leq 2$. A potential improvement in $p_{SDP}$ of the SDP solver is to explore the sparse structure of the matrices [31] or to utilize a commercial solver such as MOSEK.

# Conclusions and future research directions

From the tensor analysis of the power flow, we developed a star and linear model to achieve a scalable distributed computation. The new network model allows the direct communication between the nodal variables and the central voltages. In the model, the PL remains at unity regardless of the system size. On the other hand, the Kirchhoff's laws and voltage magnitudes are expressed in terms of nodal variables that are linear in the voltages. Therefore, the model makes it possible to keep the size of a nodal OPF small regardless of the size of the system, while the communication costs remain manageable. This new aspect of the model allows us to construct a scalable algorithm that converges to the same solution as the nonconvex OPF. We proposed the DROHS algorithm to find a local minimum using a convex surrogate function. Among the nodal OPF solutions, only near-feasible solutions (Rule 4) are selected for updating. In addition to the high quality of the solution, it also achieves computational efficiency and robustness. We tested the DROHS algorithm for the 3-, 4-, 9-, 14-, 24-, 30-, 39-, 57-, 85-, 118-, 300-, and 2,000-bus systems. The proposed algorithm achieves 1) fast and uniform convergence, 2) provable convergence, 3) the same problem formulation as the central OPF problem without ignoring any constraints, 4) guaranteed convergence to a local minimum, rather than a maximum or saddle point, that meets the first-order necessary conditions for optimality, and 5) a completely distributed algorithm (i.e., a scalable algorithm, which has never been achieved before in the literature). The challenges that the proposed algorithm faces are 1) an increased number of nodal variables that may be higher than that of the variables in the central OPF for a small system, 2) an increased number of iterations when highly connected nodes involve solutions far from the feasible regions, and 3) a prolonged wait time for nodes with low cardinalities. Therefore, the proposed algorithm is an efficient alternative to the central OPF for a large-scale network. Future research directions include the development of 1) a way to accommodate the impact of the rejected solutions in updating the x-variables if the corresponding nodes are highly connected, 2) an efficient computation to solve the nodal SDP, particularly a way to explore the sparse structure of the nodal OPF, and 3) asynchronous distributed optimization for improving the computational efficiency where the scheduling of the distributed computation is identified in terms of a knapsack problem. We also present the proof showing that the surrogate function improves at every iteration and that the iteration converges to a fixed point of the nonconvex OPF. The numerical results exhibit rapid convergence, and the convergence behavior is discussed.

# Acknowledgments

We thank Dr. Charles Van Loan who provided insight on tensor computation. We would like to express our gratitude to Dr. Robert J. Thomas and Mr. Gilbert Bindewald for their expertise on the computation for the power system analysis during this research, and we also thank 5 anonymous reviewers whose comments and suggestions helped improve and clarify this manuscript.

# Nomenclature

$\Phi_j^{M_{eq}}, \Phi_j^{M_{in}}$ eigenvectors of $M_{eq}, M_{in}$ scaled by $diag\left(\sqrt{|\lambda_{eq}|}\right)$ and $diag\left(\sqrt{|\lambda_{in}|}\right)$

$\Omega_m, \Omega_m^c$ Set and the complement of the set of Case $m$

$\alpha_j$, $\beta_j$ nodal variable associated with real and with reactive power injection, $\alpha_j^{4\times1} = \Phi_{S_j}^T v$ and $\beta_j^{4\times1} = \Phi_{\bar{S}_j}^T v$ $\gamma_{j\text{-}l,j}$, $\delta_{j\text{-}l,j}$ nodal variable with real and with reactive power flow over the line *j-l* at the side of the $j^{\text{th}}$ node, $\gamma_{j-l,j}^{4\times1} = \Phi_{F_{j-l}^{j}}^{T} v$ and $\delta_{j-l,j}^{4\times1} = \Phi_{\bar{F}_{j-l}^{j}}^{T} v$

$\omega_j$ nodal variable with $E_j$, $\omega_j^{2\times1} = \left(e_j^T + e_{j+Nb}^T\right)^T v$

*B, Br, G* sets of nodes, branches, and generators

*Br{j}, G{j}* sets of branches and of generators at node *j*

$E_j$ voltage magnitude square at Node *j*, $\left|\boldsymbol{v}_j\right|^2$

$\boldsymbol{II}_j^{eq}, \boldsymbol{II}_j^{in}$ diagonal matrix of ±1 associated with equality and inequality constraints

*I* identity matrix

***J*** matrix $[I^T\ \mathbf{j}I^T]^T$

$M_{eq}, M_{in}$ symmetric matrices with equality and inequality constraints $M_j^{eq} = \Phi_j^{M_{eq}} II_j^{eq} \Phi_j^{M_{eq}T}, M_j^{iin} = \Phi_j^{M_{in}} II_j^{in} \Phi_j^{M_{in}T}$

$\boldsymbol{Y_{bus}}$, $\boldsymbol{Y_{br}}$ nodal and branch admittance matrices

*Nb, Nl, Ng* number of nodes, lines, and generators at the system of interest

$N_{iter}$ number of iterations until convergence occurs

$\text{cap}_l$ thermal limit of flow over *l*

$\boldsymbol{d}_j$ real and reactive power demand at *j*, $\boldsymbol{d}_j = d_j + \mathbf{j}\bar{d}_j$

$e_j$ $j^{\text{th}}$ column vector in the identity matrix

$\boldsymbol{f}_{j-l}^{j}$ power flow over *j-l* at the $j^{\text{th}}$ node, $\boldsymbol{f}_{j-l}^{j} = f_{j-l}^{j} + \mathbf{j}\bar{f}_{j-l}^{j}$

$\boldsymbol{g}_j$ generation at the $j^{\text{th}}$ node, $\boldsymbol{g}_j = g_j + \mathbf{j}\bar{g}_j$

$\mathbf{j}$ $\sqrt{-1}$

$nl_j$, $ng_j$ number of lines and generators at Node *j*

$nl_{\max}$ maximum number among $nl_j$ in the system

$v_L$, $v_M$ voltages at the power and voltage channels

$\boldsymbol{v}$, $v$, $v_x$, $v_y$ complex voltage vector, voltage $v = \left[v_x^T\ v_y^T\right]^T$, real and imaginary part of voltage, $\boldsymbol{v} = v_x + \mathbf{j}v_y = \boldsymbol{J}^T v$ $\boldsymbol{x}, x, x, \overline{x}$ complex variable, real, and imaginary parts of $\boldsymbol{x}$, $\boldsymbol{x} = x + \mathbf{j}\overline{x}$ and $x = \begin{bmatrix} x^T & \overline{x}^T \end{bmatrix}^T$

# Nodal variables and their cardinality

$\alpha_j^{4\times1}$ nodal variable associated with real power injection at Node $j$

$\beta_j^{4\times1}$ nodal variable associated with reactive power injection at Node $j$

$\gamma_{jk}^{4\times1}$ nodal variable associated with real power flow $\check{f}_{j-k}^{j}$ over a line $j$-$k$ at the side of Node $j$

$\gamma_j^{4nl_j\times1}$ nodal variable associated with real power flow over the lines connected to Node $j$

$\delta_{jk}^{4\times1}$ nodal variable associated with reactive power flow $\bar{f}_{j-k}^{j}$ over a line $j$-$k$ at the side of Node $j$

$\delta_j^{4nl_j\times1}$ nodal variable associated with reactive power flow over the lines connected to Node $j$

$\omega_j^{2\times1}$ nodal variable associated with the voltage magnitude at Node $j$

$\check{f}_{j-k}^{j}$ scalar representing real power flow over a line $j$-$k$ at the side of Node $j$

$\bar{f}_{j-k}^{j}$ scalar representing reactive power over a line $j$-$k$ at the side of Node $j$

$p_j^{ng_j\times1}$ nodal real power generation vector at Node $j$

$q_j^{ng_j\times1}$ nodal reactive power generation vector at Node $j$

# Appendix

## Appendix A: Proof of claims

**Claim 1**: The matrices associated with power flows and power injections are all of rank 4.

**Proof**: The power flows and injections are the product between the voltages and current complex conjugate, $\boldsymbol{vi}^*$. The power at the $i^{\text{th}}$ node is: $p_i + \mathbf{j}q_i = \boldsymbol{v}_i\boldsymbol{i}_i^* = v^T\boldsymbol{S}_i v = v^T\left(S_i + \mathbf{j}\overline{S}_i\right)v$. Let the real-valued matrix sandwiched by $v$ be $M_i$. For $M_i$ has two nonzero rows at $i$ and $Nb+i$ rows, $M_i = 2e_i a_i^T + 2e_{i+Nb} b_i^T$ where $a_i$ and $b_i$ are column vectors. Because $p_i$ and $q_i$ are scalar, the real-valued can be replaced by symmetric $[M_i] = \frac{1}{2}\left(M_i + M_i^T\right)$ . The matrix $[M_i]$ is $[M_i] = e_i a_i^T + a_i e_i^T + e_{i+Nb} b_i^T + b_i e_{i+Nb}^T$ . The eigenvalue decomposition of the first pair of components leads to $e_j a_j^T + a_j e_j^T = \lambda_1 u_1 u_1^T + \lambda_2 u_2 u_2^T$ where $\lambda_1 = \frac{1}{2}\left\|e_j + a_j\right\|^2$, $u_1 = \frac{e_j + a_j}{\left\|e_j + a_j\right\|}$, $\lambda_2 = -\frac{1}{2}\left\|e_j - a_j\right\|^2$, and $u_2 = \frac{e_j - a_j}{\left\|e_j - a_j\right\|}$. The dimension of the null space of the matrices is ($2nb - 2$). According to the rank-nullity theorem, the rank of the first pair is $2nb - (2nb - 2) = 2$. Similarly, the rank of the last pair of two components also reveals a rank 2 decomposition. Because $M_i$ is a full row-rank matrix, the first and last pairs of components are independent (i.e., each pair spans the real and imaginary voltage spaces); hence, they do not overlap. Therefore, $[M_i]$ is a rank 4 matrix. The proof for the ranks of the matrices associated with the power flows follows the same steps.

**Claim 2**: The matrices associated with the squares of the voltage magnitudes have rank 2.

**Proof**: The voltage magnitude squared at Node $i$ is $E_j = v_{xi}^2 + v_{yi}^2 = v^T\left(e_i e_i^T + e_{i+Nb} e_{i+Nb}^T\right)v$. The eigenvector decomposition of the matrix sandwiched by the voltages reveals the rank-2 matrix.

# Appendix B: Proof of convergence

## Strategy of the proof

The proof comprises of three parts. The first part reveals the relationship between the local variables and the central variable at each iteration and proves that the local and central variables are bounded by the nodal updates. It also shows that the surrogate function always improves at each iteration. The second part shows that the surrogate function converges to a local minimum. The third part proves that the original nodal distributed OPF and the surrogate function share the same fixed point. This concludes the convergence of the proposed algorithm.

The nodal OPF problem was presented. The generalized OPF problem $w$ is described in Eq. (1) below:

$$\min_{x_j \in X_j} \sum_j w_j\left(x_j, y, z_j\right) = \min_{x_j \in X_j} \sum_j \left[ f_j\left(x_j\right) + g_j\left(x_j, y, z_j\right) \right] \tag{A1}$$

where $x$ is the nodal primal variable, $y = \left(v_L^T \quad v_M^T\right)^T$, $z$ represents the multipliers, $\Phi_j$ is the collection of $\Phi_j^{M_{eq}}$ and $\Phi_j^{M_{in}}$ corresponding to all the constraints for the $j^{th}$ nodal OPF, and $X_j$ is the column space of $\Phi_j$ at $j$.

Eq. (A1) implies that:

- $w_j$ is coercive.
- $f_j$ is nonconvex, where $f_j\left(x_j\right) \triangleq \min_{x_j \in X_i} \sum_i x_i^T M_{ob}^i x_i$, and has the following properties: 1) $f_j$ is a nonconvex but smooth function that is $C^1$ on an open set containing $X_j$ (defined by the column space of $\Phi_j$), and 2) $\nabla f_j$ is Lipschitz continuous on $X_j$ where each $X_j$ is nonempty, closed, and convex.
- $g_j$ is a quadratic ADMM function defined as $g\left(x_j, y, z_j\right) = \frac{\rho_j}{2}\left\| x_j - \Phi_j^T y \right\|^2 + z_j^T\left(x_j - \Phi_j^T y\right)$.

Even though the decision variables in (1) have low cardinality, the uniqueness and the existence of the solution are not guaranteed due to the nonconvex nature of the problem. To address the complexity issue, a surrogate function $h_j$ is introduced:

$$H\left(x, y, z\right) = \sum_j h_j\left(x_j, y, z_j\right) = \sum_j \left[ u_j\left(x_j\right) + g_j\left(x_j, y, z_j\right) \right] \quad \text{(2) where } u_j\left(x_j\right) = \text{SDP of } f_j\left(x_j\right) \tag{A2}$$

Because $u_j$ is the SDP relaxation of $f_j$, it is: 1) strongly convex, 2) Lipchitz continuous, 3) continuously differentiable on $X_j$. and 4) $\nabla_{x_j} f_j\left(x_j\right) = \nabla_{x_j} u_j\left(x_j\right)$.

## Part I: Boundedness of the updates of the variables

The ADMM-type distributed algorithm comprises four update processes – $x$- (Rule 3) and $y$-optimizations (Rule 5), as well as $z$- (Rule 5) and $x$-updates (Rule 6).

### Rule 3: x-optimization

For the solution $x_j^{k+1}$ found in Rule 3, where $\tilde{x}_j^{k+1} \triangleq \arg\min_{x_j \in X_j} h_j\left(x_j, y^k, z_j^k\right)$, the optimality condition yields:

$$\nabla_{x_j} u_j\left(x_j^{k+1}\right)+z_j^k+\rho_j\left(x_j^{k+1}-x_j^k\right)=0 \tag{A3}$$

or

$$\nabla_x U\left(x^{k+1}\right)+z^k+D_\rho\left(x^{k+1}-x^k\right)=0 \text{ in the vector form} \tag{A4}$$

where $D_\rho \triangleq diag\left(\rho_1 1_{nx_1}^T \cdots \rho_{Nb} 1_{nx_{Nb}}^T\right)$ and $\tilde{x}^{k+1}=\left(\tilde{x}_1^{(k+1)T} \cdots \tilde{x}_{Nb}^{(k+1)T}\right)^T$. We also define $\hat{x}^{k+1}=\left(\hat{x}_1^{(k+1)T} \cdots \hat{x}_{Nb}^{(k+1)T}\right)^T$, the ($k$+$1$)th update $x^{k+1}=\left(x_1^{(k+1)T} \cdots x_{Nb}^{(k+1)T}\right)^T$, and $z^{k+1}=\left(z_1^{(k+1)T} \cdots z_{Nb}^{(k+1)T}\right)^T$. The update in the $x$-variable is generalized as $\hat{x}_j^{k+1}-x_j^k=\left(1-\delta_{j3}\right)\Delta^k\left(x_j^{k+1}-x_j^k\right)+\delta_{j2}\Delta^k\varepsilon_j^{k+1}$, where $\delta_{jm}$ is the Kronecker delta that equals 1 for $j = m$ (i.e., $j \in \Omega_m$) and 0 otherwise. Rule 4 leads to:

$$\hat{\zeta}^{k+1}=x^k+\Phi^T\left(\Phi D_\rho \Phi^T\right)^{-1}\Phi D_\rho\left(x^{k+1}-x^k\right)_{\Omega_3^c} \tag{A5}$$

and

$$\hat{x}^{k+1}=x^k+\Delta^k\left(x^{k+1}-x^k\right)_{\Omega_3^c}+\Delta^k\left(\varepsilon^{k+1}\right)_{\Omega_2} \quad \text{where } \tilde{\varepsilon}^k=\left(\tilde{\varepsilon}_1^{kT} \cdots \tilde{\varepsilon}_{Nb}^{kT}\right)^T \text{ and } \tilde{\varepsilon}_j^{k+1} \triangleq \tilde{\zeta}_j^{k+1}-\tilde{x}_j^{k+1} \tag{A6}$$

## Rule 5: *y*-optimization and z-update process

$y^{k+1}$ satisfies the optimality condition $\nabla_y G\left(\hat{x}^{k+1}, y^{k+1}, z^k\right)=-\Phi D_\rho \hat{x}^{k+1}+\Phi D_\rho \Phi^T y^{k+1}-\Phi z^k=0$ in the $y$-optimization, where $\Phi=\left[\Phi_1 \cdots \Phi_{2nb}\right]$ and $z^k=\left(z_1^{kT} \cdots z_{nb}^{kT}\right)^T$. The solution is $y^{k+1}=\left(\Phi D_\rho \Phi^T\right)^{-1}\left(\Phi D_\rho \hat{x}^{k+1}+\Phi z^k\right)$. Using the solution and the $z$-update in Rule 5, $z^{k+1}=D_\rho \hat{x}^{k+1}-D_\rho \Phi^T y^{k+1}+z^k$, it is found that $z^{k+1}=\left[I-D_\rho \Phi^T\left(\Phi D_\rho \Phi^T\right)^{-1}\Phi\right]\left(D_\rho \hat{x}^{k+1}+z^k\right)$. Multiplication with the matrix $\Phi$ yields $\Phi z^{k+1}=0$ for all $k$.

$$\Delta z^{k+1}=z^{k+1}-z^k=\left[I-D_\rho \Phi^T\left(\Phi D_\rho \Phi^T\right)^{-1}\Phi\right]D_\rho \hat{x}^{k+1} \tag{A7}$$

$$y^{k+1}=\left(\Phi D_\rho \Phi^T\right)^{-1}\Phi D_\rho \hat{x}^{k+1} \tag{A8}$$

Even though the cardinality of $y$ increases with the system size, its computation in (A8) is the linear combination of nodal updates. Note that $\left(\Phi D_\rho \Phi^T\right)^{-1}\Phi D_\rho$ does not change with the iteration and that with each column in $\left(\Phi D_\rho \Phi^T\right)^{-1} \Phi D_\rho$, the nodal variable $\hat{x}^{k+1}$ finds a corresponding $y$-variable, i.e., $y^{k+1} = \sum_j y_j^{k+1}$ where $y_j^{k+1} = \left(\Phi D_\rho \Phi^T\right)^{-1} \Phi D_\rho e_j \hat{x}_j^{k+1}$. Similarly, $z^{k+1} = z^k + \Delta z^{k+1} = z^k + \sum_j \Delta \hat{z}_j^{k+1}$ where $\Delta \hat{z}_j^{k+1} = \left[ I - D_\rho \Phi^T \left(\Phi D_\rho \Phi^T\right)^{-1} \Phi \right] D_\rho e_j \hat{x}_j^{k+1}$. It is possible to compute $y_j^{k+1}$ and $\Delta \hat{z}_j^{k+1}$ at each node without an information exchange with any other nodes. While the $y$-optimization and the $z$-update process involve central information exchanges, the computations are linear updates of nodal $y_j^{k+1}$ and $\Delta \hat{z}_j^{k+1}$.

## Rule 6: $x$-update

According to Rule 6, $x$ is the projection of $\hat{x}^{k+1}$ on $X$, i.e., $x^{k+1} = \Phi^T \left(\Phi D_\rho \Phi^T\right)^{-1} \Phi D_\rho \hat{x}^{k+1}$, which yields:

$$\Phi D_\rho x^{k+1} = \Phi D_\rho \hat{x}^{k+1} \tag{A9}$$

The projection on the real space of Φ is done to make the $\hat{x}^{k+1}$ variable internally consistent with the global variable $y$, $x^{k+1} = \Phi^T y^{k+1}$. Since $x^k$ is strictly in the column space of $\Phi$, $x^k = \Phi^T \left(\Phi D_\rho \Phi^T\right)^{-1} \Phi D_\rho x^k$, i.e., the linear projection of a vector onto its full space is the vector itself, and

$$x^{k+1} - x^k = \Phi^T \left(\Phi D_\rho \Phi^T\right)^{-1} \Phi D_\rho \left(\hat{x}^{k+1} - x^k\right) \tag{A10}$$

$$y^{k+1} - y^k = \left(\Phi D_\rho \Phi^T\right)^{-1} \Phi D_\rho \left(x^{k+1} - x^k\right) \tag{A11}$$

$$z^{k+1} - z^k = \left[ I - D_\rho \Phi^T \left(\Phi D_\rho \Phi^T\right)^{-1} \Phi \right] D_\rho \left(\hat{x}^{k+1} - x^k\right) \tag{A12}$$

Eq. (A10) indicates that a linear mapping exists between $x^{k+1} - x^k$ and $\hat{x}^{k+1} - x^k$. Eq. (A12) implies that $z^{k+1} - z^k \in \mathrm{Imag}\left[ I - D_\rho \Phi^T \left(\Phi D_\rho \Phi^T\right)^{-1} \Phi \right]$. Note that $\left[ I - D_\rho^{1/2} \Phi^T \left(\Phi D_\rho \Phi^T\right)^{-1} \Phi D_\rho^{1/2} \right] D_\rho^{1/2} \Phi^T = 0$, i.e., $\Phi D_\rho^{1/2}$ exists in the null space of $\left[ I - D_\rho^{1/2} \Phi^T \left(\Phi D_\rho \Phi^T\right)^{-1} \Phi D_\rho^{1/2} \right]$. According to Rule 6, i.e., $x^k = \Phi^T y^k$ for all $k$, $\Delta x^{k+1} \left(\triangleq x^{k+1} - x^k\right)$ is also in the column space of Φ. Then, $\left[ I - D_\rho^{1/2} \Phi^T \left(\Phi D_\rho \Phi^T\right)^{-1} \Phi D_\rho^{1/2} \right] D_\rho^{1/2} \Delta x^{k+1} = 0$. From Eq. (A7), we arrive at $\left(\Delta z^{k+1}\right)^T D_\rho^{1/2} \Delta x^{k+1} = 0$, meaning that the update of the primal variables and the multiplier update are perpendicular. Hence:

$D_\rho^{1/2}\Delta x^{k+1} = D_\rho^{1/2}\Phi^T\left(\Phi D_\rho \Phi^T\right)^{-1}\Phi D_\rho^{1/2} D_\rho^{1/2}\Delta x^{k+1}$. The singular value decomposition of $\Phi D_\rho^{1/2}$, i.e., $\Phi D_\rho^{1/2} = U\Sigma V^T$, yields: $D_\rho^{1/2}\Delta x^{k+1} = V\begin{bmatrix}\Sigma_R \\ 0\end{bmatrix}U^T\left(U\Sigma_R^2 U^T\right)^{-1}U\begin{bmatrix}\Sigma_R & 0\end{bmatrix}V^T D_\rho^{1/2}\Delta x^{k+1}$ which is further simplified $D_\rho^{1/2}\Delta x^{k+1} = V_R V_R^T D_\rho^{1/2}\Delta x^{k+1}$. Therefore,

$$V^T D_\rho^{1/2}\left(x^{k+1} - x^k\right) = V_R^T D_\rho^{1/2}\left(x^{k+1} - x^k\right) \tag{A13}$$

Eq. (A13) means that $D_\rho^{1/2}\left(x^{k+1} - x^k\right)$ stays strictly in the real space of $V$ or of $\Phi D_\rho^{1/2}$.

On the other hand, from Eq (A10), it can be seen that $\Phi D_\rho^{1/2}\left[D_\rho^{1/2}\left(x^{k+1} - \hat{x}^{k+1}\right)\right] = 0$, i.e., $D_\rho^{1/2}\left(x^{k+1} - \hat{x}^{k+1}\right)$ stays strictly in the null space of $V$ or of $\Phi D_\rho^{1/2}$. Using Eq. (A13), it can be seen that $D_\rho^{1/2}\left(x^{k+1} - x^k\right)$ and $D_\rho^{1/2}\left(x^{k+1} - \hat{x}^{k+1}\right)$ are perpendicular. According to (A6),

$$\left\|\hat{x}^{k+1} - x^k\right\| = \Delta^k\left\|\left(x^{k+1} - x^k\right)_{\Omega_3^c} + \left(\varepsilon^{k+1}\right)_{\Omega_2}\right\| \le \Delta^k\left\|x^{k+1} - x^k\right\|_{\Omega_3^c} + \Delta^k\left\|\varepsilon^{k+1}\right\|_{\Omega_2} \le \Delta^k\left(1 + \tau_{\max}^k\right)\left\|x^{k+1} - x^k\right\|_{\Omega_3^c} \tag{A14}$$

Eqs. (A10), (A11), (A12), and (A14) yield:

$$\left\|x^{k+1} - x^k\right\| = \left\|D_\rho^{-1/2} V\Sigma^T\Sigma_R^{-2}\Sigma V^T D_\rho^{1/2}\left(\hat{x}^{k+1} - x^k\right)\right\| \le \Delta^k\left(1 + \tau_{\max}^k\right)\sqrt{\frac{\rho_{\max}}{\rho_{\min}}}\left\|x^{k+1} - x^k\right\|_{\Omega_3^c} \tag{A15}$$

$$\left\|y^{k+1} - y^k\right\| = \left\|U\begin{bmatrix}\Sigma_R^{-1} & 0\end{bmatrix}V^T D_\rho^{1/2}\left(x^{k+1} - x^k\right)\right\| \le \frac{\Delta^k\left(1 + \tau_{\max}^k\right)\sqrt{\rho_{\max}}}{\left|\sigma_{\Phi D_\rho^{1/2}}^{\min}\right|}\left\|x^{k+1} - x^k\right\|_{\Omega_3^c} \tag{A16}$$

$$\left\|z^{k+1} - z^k\right\| = \left\|D_\rho^{1/2} V\left(I - \Sigma^T\Sigma_R^{-2}\Sigma\right)V^T D_\rho^{1/2}\left(\hat{x}^{k+1} - x^k\right)\right\| \le \Delta^k\left(1 + \tau_{\max}^k\right)\rho_{\max}\left\|x^{k+1} - x^k\right\|_{\Omega_3^c} \tag{A17}$$

Inequalities (A15), (A16), and (A17) show that the nodal updates are all bound by finite multiples of $\Delta^k\left\|x^{k+1} - x^k\right\|_{\Omega_3^c}$. Note that $\lim_{k\to\infty}\Delta^k = 0$ for the series $\Delta^{k+1} = \Delta^k - a\left(\Delta^k\right)^2$.

# Part II: Convergence of the surrogate function

In this section, we will show that in the consecutive updates in $x$- (Rule 3) and $y$-optimizations (Rule 5), $z$- (Rule 5), and $x$-updates (Rule 6), the convex surrogate function converges uniformly.

## The convergence of the surrogate function *H* to *W*

$G\left(x^{k+1},y,^{k+1},z^{k+1}\right)-G\left(x^k,y^k,z^k\right)=0$. For $G\left(x^k,y^k,z^k\right)\triangleq\frac{1}{2}\left\|D_\rho^{1/2}\left(x^k-\Phi^T y^k\right)\right\|^2+z^{kT}\left(x^k-\Phi^T y^k\right)$, the change in $H$ becomes:

$$H\left(x^{k+1},y,^{k+1},z^{k+1}\right)-H\left(x^k,y^k,z^k\right)=U\left(x^{k+1}\right)-U\left(x^k\right) \tag{A18}$$

The strong convexity of $U$ in $[x^k, x^{k+1}]$ yields:

$$\begin{aligned} U\left(x^{k+1}\right)&\le U\left(x^k\right)+\Delta^k\left[U\left(\hat{\zeta}^{k+1}+\hat{\varepsilon}^{k+1}\right)-U\left(x^k\right)\right]-\frac{\Delta^k\left(1-\Delta^k\right)c_U}{2}\left\|\hat{\zeta}^{k+1}+\hat{\varepsilon}^{k+1}-x^k\right\|^2\\ &=U\left(x^k\right)+\Delta^k\left[U\left(\hat{\zeta}^{k+1}\right)-U\left(x^k\right)\right]-\frac{\Delta^k\left(1-\Delta^k\right)c_U}{2}\left\|\hat{\zeta}^{k+1}-x^k\right\|^2+\vartheta\left(\Delta^k\tau_{\max}^k\right)\end{aligned} \tag{A19}$$

and that in [ $x^k$, $\hat{\zeta}^{k+1}$] yields:

$$U\left(\hat{\zeta}^{k+1}\right)-U\left(x^k\right)\le\left[\nabla_x U\left(\hat{\zeta}^{k+1}\right)\right]^T\left(\hat{\zeta}^{k+1}-x^k\right)-\frac{c_U}{2}\left\|\hat{\zeta}^{k+1}-x^k\right\|^2 \tag{A20}$$

Inequalities (A19) and (A20) yield:

$$U\left(x^{k+1}\right)\le U\left(x^k\right)+\Delta^k\left[\nabla_x U\left(\hat{\zeta}^{k+1}\right)\right]^T\left(\hat{\zeta}^{k+1}-x^k\right)-\Delta^k c_U\left(1-\frac{\Delta^k}{2}\right)\left\|\hat{\zeta}^{k+1}-x^k\right\|^2+\vartheta\left(\tau_{\max}^k\right) \tag{A21}$$

The Lipschitz continuity of $\nabla_x U\left(x\right)$ yields $\left[\nabla_x U\left(\hat{\zeta}^{k+1}\right)-\nabla_x U\left(x^k\right)\right]^T\left(\hat{\zeta}^{k+1}-x^k\right)\le L_{\nabla U}\left\|\hat{\zeta}^{k+1}-x^k\right\|^2$, and Inequality (A21) leads to:

$$U\left(x^{k+1}\right)\le U\left(x^k\right)+\Delta^k\left[\nabla_x U\left(x^k\right)\right]^T\left(\hat{\zeta}^{k+1}-x^k\right)-c_U\Delta^k\left(1-\frac{\Delta^k}{2}-\frac{L_{\nabla U}}{c_U}\right)\left\|\hat{\zeta}^{k+1}-x^k\right\|^2+\vartheta\left(\tau_{\max}^k\right) \tag{A22}$$

Using Eq. (A6), the second term in Inequality (A22) is written as follows:

$$\left[\nabla_x U\left(x^k\right)\right]^T\left(\hat{\zeta}^{k+1}-x^k\right)=\left[\nabla_x U\left(x^k\right)\right]^T\Phi^T\left(\Phi D_\rho\Phi^T\right)^{-1}\Phi D_\rho\left(x^{k+1}-x^k\right)_{\Omega_3^c} \tag{A23}$$

Using Eq. (A3) and $\Phi z^k = 0$, we arrive at:

$$\begin{aligned} &\left[\nabla_x U\left(x^{k+1}\right)\right]^T\Phi^T\left(\Phi D_\rho\Phi^T\right)^{-1}\Phi D_\rho\left(x^{k+1}-x^k\right)_{\Omega_3^c}=-\left(x^{k+1}-x^{k+1}\right)^T D_\rho\Phi^T\left(\Phi D_\rho\Phi^T\right)^{-1}\Phi D_\rho\left(x^{k+1}-x^k\right)_{\Omega_3^c}\\ &=-\left(x^{k+1}-x^k\right)^T D_\rho\Phi^T\left(\Phi D_\rho\Phi^T\right)^{-1}\Phi D_\rho\left(x^{k+1}-x^k\right)_{\Omega_3^c}+\left(x^{k+1}-x^k\right)^T D_\rho\Phi^T\left(\Phi D_\rho\Phi^T\right)^{-1}\Phi D_\rho\left(x^{k+1}-x^k\right)_{\Omega_3^c}\\ &=-\left\|V_R^T D_\rho^{1/2}\left(x^{k+1}-x^k\right)\right\|_{\Omega_3^c}^2+\left(y^{k+1}-y^k\right)^T\Phi D_\rho\left(x^{k+1}-x^k\right)_{\Omega_3^c}=-\left\|V_R^T D_\rho^{1/2}\left(x^{k+1}-x^k\right)\right\|_{\Omega_3^c}^2+\left(x^{k+1}-x^k\right)^T D_\rho\left(x^{k+1}-x^k\right)_{\Omega_3^c}\\ &=-\left\|V_R^T D_\rho^{1/2}\left(x^{k+1}-x^k\right)\right\|_{\Omega_3^c}^2+\left\|D_\rho^{1/2}\left(x^{k+1}-x^k\right)\right\|^2\end{aligned} \tag{A24}$$

Ineq. (A15) and Eq. (A24) yield:

$$\left[\nabla_x U\left(x^{k+1}\right)\right]^T \Phi^T\left(\Phi D_\rho \Phi^T\right)^{-1} \Phi D_\rho\left(x^{k+1}-x^k\right)_{\Omega_3^c} \leq -\rho_{\min}\left\|x^{k+1}-x^k\right\|_{\Omega_3^c}^2+\left(\Delta^k\right)^2 \rho_{\max}\left\|x^{k+1}-x^k\right\|_{\Omega_3^c}^2 \tag{A25}$$

By the definition of $\hat{\zeta}^{k+1}$, it follows that:

$$\left\|\hat{\zeta}^{k+1}-x^k\right\|=\left\|\Phi^T\left(\Phi D_\rho \Phi^T\right)^{-1} \Phi D_\rho\left(x^{k+1}-x^k\right)\right\|_{\Omega_3^c} \leq \sqrt{\frac{\rho_{\max}}{\rho_{\min}}}\left\|x^{k+1}-x^k\right\|_{\Omega_3^c} \tag{A26}$$

Ineq. (A22) becomes:

$$U\left(x^{k+1}\right) \leq U\left(x^k\right)-\Delta^k\left[\rho_{\min}+\frac{\rho_{\max}}{\rho_{\min}} c_U-\left(\frac{\Delta^k}{2} c_U+L_{\nabla U}\right) \frac{\rho_{\max}}{\rho_{\min}}-\left(\Delta^k\right)^2 \rho_{\max}\right]\left\|x^{k+1}-x^k\right\|_{\Omega_3^c}^2+\vartheta\left(\tau_{\max}^k\right) \tag{A27}$$

For a given $\Delta^k$ in (0, $\Delta_{max}$] and a given ratio $\rho_{min}/\rho_{max}$, there exists $\rho_{\min}^{cr}$ such that $\rho_{\min}^{cr}=\max\left\{0, \frac{L_{\nabla U}-\left(1-\frac{\Delta_{\max}}{2}\right) c_U}{\frac{\rho_{\min}}{\rho_{\max}}-\Delta_{\max}^2}\right\}$.

From Rule 5 in the proposed algorithm, $\Delta^{k+1}=-a\left(\Delta^k-\frac{1}{2a}\right)^2+\frac{1}{4a}$, which leads to $\Delta^{k+1} \leq \frac{1}{4a}$, and the equality holds when $\Delta^k=\frac{1}{2a}$. As $k$ increases, Rule 5 ensures that $\Delta^k$ decreases monotonously. If $a$ is set at less than ¼, $\rho_{\min}^{cr}$ is strictly positive. Since the ADMM-type constraint function $G$ does not change at the iteration, we have the following:

$$H\left(x^{k+1}, y^{k+1}, z^{k+1}\right) \leq H\left(x^k, y^k, z^k\right)-\Delta^k\left[\rho_{\min}+\frac{\rho_{\max}}{\rho_{\min}} c_U-\left(\frac{\Delta^k}{2} c_U+L_{\nabla U}\right) \frac{\rho_{\max}}{\rho_{\min}}-\left(\Delta^k\right)^2 \rho_{\max}\right]\left\|x^{k+1}-x^k\right\|_{\Omega_3^c}^2+\vartheta\left(\tau_{\max}^k\right) \tag{A28}$$

For a sufficiently large $\rho > \rho_{\min}^{cr}$ (> 0), the objective function $H$'s surrogate function always decreases with each iteration. Because $H$ is coercive and always improves at iteration, the surrogate function converges.

## Part III: Convergence of the proposed algorithm

In this section, we present the proof of the convergence of the proposed algorithm by showing: 1) the series of solutions to the surrogate function found through iteration is Lipchitz continuous; 2) the surrogate function and the nonconvex OPF share a fixed point; and 3) the proposed algorithm converges to the shared fixed point, which is a solution to the nonconvex OPF.

## Lipschitz continuity of $\chi_j^{k+1}\left(\chi_j^k\right)$ and the fixed point

Let $\chi_j$ be the aggregated nodal variable ($x_j$, $y$, $z_j$). For an initial point $\chi^0$, $x_j\left(\chi^0\right) \quad \underset{x_j}{\arg\min}\, h_j\left(x_j \middle| \chi^0\right)$ is uniquely defined. All $\chi_j^0$ share the global voltage $y$ only and are otherwise independent, which makes that $\tilde{x}_j\left(\chi^0\right)=\tilde{x}_j\left(\chi_j^0\right)\triangleq \underset{x_j}{\arg\min}\, h_j\left(x_j \middle| \chi_j^0\right)$. Since $x_j$, $y$, and $z_j$ are all determined by $x_j$, $\chi_j$ is the unique solution to $\tilde{\chi}_j\left(\chi_j^0\right)\triangleq \underset{\chi_j}{\arg\min}\, h_j\left(\chi_j \middle| \chi_j^0\right)$. For the unique solution $\chi_j$ and a convex function $h_j$, where $h_j\left[\chi_j\left(\chi_j^0\right)\right]\le h_j\left(\chi_j^0\right)$, the optimality condition provides the following result:

$$\left\{\nabla_{\chi_j} h_j\left[\chi_j\left(\chi_j^0\right)\right]\right\}^T\left[\chi_j^0-\chi_j\left(\chi_j^0\right)\right]\ge 0 \tag{A29}$$

By the definition of $h_j$ in Eq. (A2), the optimality condition is:

$$\left\{\nabla_{\chi_j} h_j\left[\chi_j\left(\chi_j^0\right)\right]\right\}^T\left[\chi_j^0-\chi_j\left(\chi_j^0\right)\right]=\left\{\nabla_{\chi_j} u_j\left[\chi_j\left(\chi_j^0\right)\right]+\nabla_{\chi_j} g_j\left[\chi_j\left(\chi_j^0\right)\right]\right\}^T\left[\chi_j^0-\chi_j\left(\chi_j^0\right)\right]\ge 0 \tag{A30}$$

The convexity of $g_j$ yields:

$$g_j\left(\chi_j^0\right)\ge g_j\left[\chi_j\left(\chi_j^0\right)\right]+\left\{\nabla_{\chi_j} g_j\left[\chi_j\left(\chi_j^0\right)\right]\right\}^T\left[\chi_j^0-\chi_j\left(\chi_j^0\right)\right]+\frac{c_g}{2}\left\|\chi_j^0-\chi_j\left(\chi_j^0\right)\right\|^2 \tag{A31}$$

Ineq. (A30) becomes

$$\left\{\nabla_{\chi_j} u_j\left[\chi_j\left(\chi_j^0\right)\right]\right\}^T\left[\chi_j^0-\chi_j\left(\chi_j^0\right)\right]+g_j\left(\chi_j^0\right)-g_j\left[\chi_j\left(\chi_j^0\right)\right]-\frac{c_g}{2}\left\|\chi_j^0-\chi_j\left(\chi_j^0\right)\right\|^2\ge\left\{\nabla_{\chi_j} h_j\left[\chi_j\left(\chi_j^0\right)\right]\right\}^T\left[\chi_j^0-\chi_j\left(\chi_j^0\right)\right]\ge 0$$

$$\rightarrow\left\{\nabla_{\chi_j} u_j\left[\chi_j\left(\chi_j^0\right)\right]\right\}^T\left[\chi_j^0-\chi_j\left(\chi_j^0\right)\right]+g_j\left(\chi_j^0\right)-g_j\left[\chi_j\left(\chi_j^0\right)\right]\ge\frac{c_g}{2}\left\|\chi_j^0-\chi_j\left(\chi_j^0\right)\right\|^2 \tag{A32}$$

Note that $\nabla_{\chi_j} u_j\left[\chi_j\left(\chi_j^0\right)\right]=\nabla_{\chi_j} f_j\left[\chi_j\left(\chi_j^0\right)\right]$, which implies that $\chi_j\left(\chi_j^0\right)$ satisfies the optimality condition for the exact and nonconvex OPF as well. Consider two starting points $\chi^I$ and $\chi^{II}$:

$$\begin{aligned}&\left\{\nabla_{\chi_j} u_j\left[\chi_j\left(\chi_j^I\right)\right]\right\}^T\left[\chi_j^0-\chi_j\left(\chi_j^I\right)\right]+g_j\left(\chi_j^0\right)-g_j\left[\chi_j\left(\chi_j^I\right)\right]\ge 0\\&\left\{\nabla_{\chi_j} u_j\left[\chi_j\left(\chi_j^{II}\right)\right]\right\}^T\left[\chi_j^0-\chi_j\left(\chi_j^{II}\right)\right]+g_j\left(\chi_j^0\right)-g_j\left[\chi_j\left(\chi_j^{II}\right)\right]\ge 0\end{aligned} \tag{A33}$$

The inequalities hold for any feasible $\chi_j^0$. Set $\chi_j^0$ to $\chi_j\left(\chi_j^{II}\right)$ (top) and to $\chi_j\left(\chi_j^I\right)$ (bottom), and add inequalities:

$$\left\{\nabla_{\chi_j} u_j\left[\chi_j\left(\chi_j^{II}\right)\right]-\nabla_{\chi_j} u_j\left[\chi_j\left(\chi_j^I\right)\right]\right\}^T\left[\chi_j\left(\chi_j^{II}\right)-\chi_j\left(\chi_j^I\right)\right]\le 0 \tag{A34}$$

The function $u_j^{\nabla}\left(\chi_j\right) \triangleq \nabla_{\chi_j} u_j\left[\tilde{\chi}_j\left(\chi_j\right)\right]$ is Lipschitz continuous because $u_j$ is continuously differentiable, i.e., $\left\|\nabla_{\chi_j} u_j\left[\chi_j\left(\chi_j^{II}\right)\right]-\nabla_{\chi_j} u_j\left[\chi_j\left(\chi_j^{I}\right)\right]\right\|=\left\|u_j^{\nabla}\left(\chi_j^{II}\right)-u_j^{\nabla}\left(\chi_j^{I}\right)\right\| \leq L_{u_j^{\nabla}}\left\|\chi_j^{II}-\chi_j^{I}\right\|$. The upper bound of Ineq. (A34) is:

$$\begin{aligned}&\left\{\nabla_{\chi_j} u_j\left[\chi_j\left(\chi_j^{II}\right)\right]-\nabla_{\chi_j} u_j\left[\chi_j\left(\chi_j^{I}\right)\right]\right\}^T\left[\chi_j\left(\chi_j^{II}\right)-\chi_j\left(\chi_j^{I}\right)\right] \\ &\quad \leq\left\|\nabla_{\chi_j} u_j\left[\chi_j\left(\chi_j^{II}\right)\right]-\nabla_{\chi_j} u_j\left[\chi_j\left(\chi_j^{I}\right)\right]\right\|\left\|\chi_j\left(\chi_j^{II}\right)-\chi_j\left(\chi_j^{I}\right)\right\| \leq L_{u_j^{\nabla}}\left\|\chi_j^{II}-\chi_j^{I}\right\|\left\|\chi_j\left(\chi_j^{II}\right)-\chi_j\left(\chi_j^{I}\right)\right\|\end{aligned} \tag{A35}$$

On the other hand, the strong convexity of $u_j\left(\chi_j\right)$ yields

$$\left\{\nabla_{\chi_j} u_j\left[\chi_j\left(\chi_j^{II}\right)\right]-\nabla_{\chi_j} u_j\left[\chi_j\left(\chi_j^{I}\right)\right]\right\}^T\left[\chi_j\left(\chi_j^{II}\right)-\chi_j\left(\chi_j^{I}\right)\right] \geq c_u\left\|\chi_j\left(\chi_j^{II}\right)-\chi_j\left(\chi_j^{I}\right)\right\|^2 \tag{A36}$$

Combining Ineq. (A35) and (A36) yields:

$$\left\|\chi_j\left(\chi_j^{II}\right)-\chi_j\left(\chi_j^{I}\right)\right\| \leq \frac{L_{u_j^{\nabla}}}{c_u}\left\|\chi_j^{II}-\chi_j^{I}\right\| \tag{A37}$$

The inequality proves the Lipschitz continuity of $\chi_j^{k+1}\left(\chi_j^k\right)$, i.e., the neighborhood of initial points stays close in the mapping of the function $\chi_j^{k+1}\left(\chi_j^k\right)$. Ineq. (A37) holds for the fixed point $\chi_j^*$, where $\chi_j^*\left(\chi_j^*\right)=\chi_j^*$: $\left[\nabla_{\chi_j} u_j\left(\chi_j^*\right)\right]^T\left(\chi_j^0-\chi_j^*\right)+g_j\left(\chi_j^0\right)-g_j\left(\chi_j^*\right) \geq 0$. For some $g_j^{\nabla} \in \partial g_j\left(\chi^*\right)$, the optimality condition reduces to: $\left[\nabla_{\chi_j} u_j\left(\chi_j^*\right)+g_j^{\nabla}\right]^T\left(\chi_j^0-\chi_j^*\right) \geq 0$. Summing over all $j$ yields $\left[\nabla_\chi U\left(\chi^*\right)+G_j^{\nabla}\right]^T\left(\chi^0-\chi^*\right) \geq 0$. Note that $\nabla_\chi U(\chi)=\nabla_\chi F(\chi)$, which implies that $\chi^*$ is a stationary point of the nonconvex OPF problem as well. $\chi^*$ is the unique solution to the SDP with initial point of $\chi^*$, and $\chi^*$ satisfies the optimality conditions. Therefore, the nonconvex OPF and its relaxed SDP share a fixed point.

## Convergence of the nonconvex OPF

For the $k^{\text{th}}$ iteration and an initial point $\chi_j^0=\chi_j^k$, Ineq. (A36) leads to: $\left\{\nabla_{\chi_j} u_j\left[\chi_j\left(\chi_j^k\right)\right]\right\}^T\left[\chi_j^k-\chi_j\left(\chi_j^k\right)\right]+g_j\left(\chi_j^k\right)-g_j\left[\chi_j\left(\chi_j^k\right)\right] \geq \frac{c_g}{2}\left\|\chi_j^k-\chi_j\left(\chi_j^k\right)\right\|^2$. By $\nabla_{\chi_j} f_j\left[\chi_j\left(\chi_j^k\right)\right]=\nabla_{\chi_j} u_j\left[\chi_j\left(\chi_j^k\right)\right]$, $\left\{\nabla_{\chi_j} f_j\left[\chi_j\left(\chi_j^k\right)\right]\right\}^T\left[\chi_j^k-\chi_j\left(\chi_j^k\right)\right]+g_j\left(\chi_j^k\right)-g_j\left[\chi_j\left(\chi_j^k\right)\right] \geq \frac{c_g}{2}\left\|\chi_j^k-\chi_j\left(\chi_j^k\right)\right\|^2$. Summing over all $j$ gives:

$$\left\{\nabla_\chi F\left[\chi\left(\chi^k\right)\right]\right\}^T\left[\chi^k-\chi\left(\chi^k\right)\right]+\sum_j\left\{g_j\left(\chi_j^k\right)-g_j\left[\chi_j\left(\chi_j^k\right)\right]\right\}\geq\frac{c_g^{\min}}{2}\left\|\chi^k-\chi\left(\chi^k\right)\right\|^2 \tag{A38}$$

Since $x_j$, $y$, and $z_j$ are all determined from $x_j$, $\chi_j$ is well defined in terms of $x_j$. At given $y^k$ and $z^k$,

$$\left\{\nabla_\chi F\left[\chi\left(x^{k+1}\right)\right]\right\}^T\left[\chi\left(x^k\right)-\chi\left(x^{k+1}\right)\right]+\sum_j\left\{g_j\left[\chi_j\left(x^k\right)\right]-g_j\left[\chi_j\left(x^{k+1}\right)\right]\right\}\geq\frac{c_g^{\min}}{2}\left\|\chi\left(x^k\right)-\chi\left(x^{k+1}\right)\right\|^2 \tag{A39}$$

Note that $F$ is a function of $x$ only and that $\chi_j$ is ($x_j$, $y$, $z_j$). Ineq. (A39) becomes:

$$\left[\nabla_x F\left(x^{k+1}\right)\right]^T\left(x^k-x^{k+1}\right)+\sum_j\left[g_j\left(x^k,y^k,z^k\right)-g_j\left(x^{k+1},y^{k+1},z^{k+1}\right)\right]\geq\frac{c_g^{\min}}{2}\left\|\begin{pmatrix}x^k-x^{k+1}\\ y^k-y^{k+1}\\ z^k-z^{k+1}\end{pmatrix}\right\|^2 \tag{A40}$$

Using Ineqs. (A15), (A16), (A17), and $\left\|I-D_\rho^{1/2}\Phi^T\left(\Phi D_\rho\Phi^T\right)^{-1}\Phi D_\rho^{1/2}\right\|=1$, it can be seen that:

$$\left\|\chi\left(x^k\right)-\chi\left(x^{k+1}\right)\right\|^2\leq\left[\left\|D_\rho^{-1}\right\|\left\|D_\rho\right\|+\left\|\left(\Phi D_\rho\Phi^T\right)^{-1}\Phi D_\rho\right\|^2+\left\|D_\rho\right\|\right]\left\|\hat{x}^{k+1}-x^k\right\|^2 \tag{A41}$$

Using $\hat{x}^{k+1}-x^k=\Delta^k\left(I-\Lambda_3\right)\left(x^{k+1}-x^k\right)+\Delta^k\Lambda_2 b^{k+1}$,

$$\left\|\hat{x}^{k+1}-x^k\right\|\leq\Delta^k\left\|x^{k+1}-x^k\right\|_{\Omega_3^c}+\Delta^k\tau_{\max}^k\left\|x^{k+1}-x^k\right\|_{\Omega_2}\leq\Delta^k\left(1+\tau_{\max}^k\right)\left\|x^{k+1}-x^k\right\|_{\Omega_3^c} \tag{A42}$$

Ineq. (A40) becomes:

$$\left[\nabla_x F\left(x^{k+1}\right)\right]^T\left(x^k-x^{k+1}\right)+\sum_j\left[g_j\left(x^k,y^k,z^k\right)-g_j\left(x^{k+1},y^{k+1},z^{k+1}\right)\right]\geq\frac{c_g^{\min}}{2}C_W^{k+1}\left(\Delta^k\right)^2\left(1+\tau_{\max}^k\right)^2\left\|x^{k+1}-x^k\right\|_{\Omega_3^c}^2 \tag{A43}$$

where $C_W^{k+1}=\left\|D_\rho^{-1}\right\|\left\|D_\rho\right\|+\left\|\left(\Phi D_\rho\Phi^T\right)^{-1}\Phi D_\rho\right\|^2+\left\|D_\rho\right\|$.

Since $U$ is strongly convex, $\left[\nabla_x U\left(x^k\right)\right]^T\times\left(x^k-x^{k+1}\right)-\frac{c_U}{2}\left\|x^{k+1}-x^k\right\|^2\geq\left[\nabla_x U\left(x^{k+1}\right)\right]^T\left(x^k-x^{k+1}\right)$. Ineq. (A43) becomes:

$$\left[\nabla_x F\left(x^k\right)\right]^T\left(x^{k+1}-x^k\right)+\sum_j\left[g_j\left(x^{k+1},y^{k+1},z^{k+1}\right)-g_j\left(x^k,y^k,z^k\right)\right]\leq-\frac{c_g^{\min}C_W^{k+1}+c_U\left\|D_\rho^{-1}\right\|\left\|D_\rho\right\|}{2}\left(\Delta^k\right)^2\left(1+\tau_{\max}^k\right)^2\left\|x^{k+1}-x^k\right\|_{\Omega_3^c}^2 \tag{A44}$$

For the $k^{\text{th}}$ iteration, the Descent Lemma results in [32]:

$$\begin{aligned}F\left(x^{k+1}\right)&\leq F\left(x^k\right)+\nabla_x F\left(x^k\right)^T\left(x^{k+1}-x^k\right)+\frac{L_F}{2}\left\|x^{k+1}-x^k\right\|^2\\&\leq F\left(x^k\right)+\nabla_x F\left(x^k\right)^T\left(x^{k+1}-x^k\right)+\frac{L_F}{2}\left(\Delta^k\right)^2\left(1+\tau_{\max}\right)^2\left\|D_\rho^{-1}\right\|\left\|D_\rho\right\|\left\|x^{k+1}-x^k\right\|_{\Omega_3^c}^2\end{aligned} \tag{A45}$$

Combining (A44), and (A45) yields:

$$
\begin{aligned}
W\left(\chi^{k+1}\right) &= F\left(x^{k+1}\right)+\sum_j g_j\left(\chi^{k+1}\right) \\
&\le F\left(x^k\right)+\sum_j g_j\left(\chi^k\right)+\nabla_x F\left(x^k\right)^T\left(x^{k+1}-x^k\right)+\sum_j g_j\left(\chi^{k+1}\right)-\sum_j g_j\left(\chi^k\right)+\frac{L_F}{2}\left(\Delta^k\right)^2\left(1+\tau_{\max}^k\right)^2\left\|D_\rho^{-1}\right\|\left\|D_\rho\right\|\left\|x^{k+1}-x^k\right\|_{\Omega_3^c}^2 \\
&\le W\left(\chi^k\right)-\eta_w^{k+1}\left(\Delta^k\right)^2\left\|x^{k+1}-x^k\right\|_{\Omega_3^c}^2+\tau_{\max}^k\left[-\eta_w^{k+1}\left(2+\tau_{\max}\right)\left(\Delta^k\right)^2\right]\left\|x^{k+1}-x^k\right\|_{\Omega_3^c}^2
\end{aligned}
$$

where $\eta_w^{k+1}=\dfrac{c_g^{\min}C_W^{k+1}+\left\|D_\rho^{-1}\right\|\left\|D_\rho\right\|\left(c_U-L_F\right)}{2}$ (A46)

The sequence $\tau_{\max}^k\left[-\eta_w^{k+1}\left(2+\tau_{\max}^k\right)\left(\Delta^k\right)^2\right]\left\|x^{k+1}-x^k\right\|_{\Omega_3^c}^2$ is bounded by $\tau_{\max}\left[-\eta_w^{k+1}\left(2+\tau_{\max}^k\right)\left(\Delta^k\right)^2\right]\left\|\varepsilon^k\right\|^2$. It is noted that $\sum_{k=1}^{\infty}\left\|\varepsilon^k\right\|^2\le\left\|x_j^{k+1}-x_j^k\right\|_{\max}^2\sum_{k=1}^{\infty}\left(\tau_{\max}^k\right)^2<\infty$ for a sequence $\tau_{\max}^k=\tau_{\max}^0/k$. Since $\sum_{k=1}^{\infty}\left\|\varepsilon^k\right\|^2<\infty$, the sum of the sequence $\tau_{\max}\left[-\eta_w^{k+1}\left(2+\tau_{\max}^k\right)\left(\Delta^k\right)^2\right]\left\|x^{k+1}-x^k\right\|_{\Omega_3^c}^2$ is finite.

Using the Lemma 3.4 in [33], either $W(\chi^k)$ divulges to -∞ or $W(\chi^k)$ converges to a finite value and $\sum_{k=1}^{\infty}\eta_w^{k+1}\left(\Delta^k\right)^2\left\|x^{k+1}-x^k\right\|_{\Omega_3^c}^2<\infty$. $W$ is coercive, which implies that $W(\chi^k)$ does not divulge. Therefore, $W(\chi^k)$ converges to a finite value, and $\sum_{k=1}^{\infty}\eta_w^{k+1}\left(\Delta^k\right)^2\left\|x^{k+1}-x^k\right\|_{\Omega_3^c}^2<\infty$. It is proven that $\lim_{k\to\infty}\left\|x^{k+1}-x^k\right\|_{\Omega_3^c}^2=0$ in [24], which means that the sequence $\{\chi^k\}$, using an inexact search, finds a local solution to the nonconvex OPF.